\documentclass[paper=a4,abstract=true,numbers=noenddot]{scrartcl}

\usepackage[english]{babel}  
\usepackage{hyperref}
\usepackage[intlimits]{amsmath}
\usepackage{amsthm}
\usepackage{amssymb} 
\usepackage{units} 
\usepackage{amsfonts,mathrsfs,bbm,mathabx}
\allowdisplaybreaks[1]
\usepackage[numbers]{natbib}
\usepackage{graphicx,subcaption}
\usepackage[space]{grffile} % + graphicx: allow figure names with spaces
\usepackage{tikz}
\usetikzlibrary{calc}
\usepackage{pgfplots}
\pgfplotsset{every axis/.append style={thick}}
\pgfplotsset{every axis legend/.append
  style={cells={anchor=west},anchor=west}}
\pgfplotscreateplotcyclelist{mylines}{%
  {red,mark=*},
  {blue,mark=square},
  {green,mark=+},
  {black,mark=o},
  {violet,mark=otimes},
  {brown,mark=triangle},
  {cyan,mark=diamond},
  {orange,mark=10-pointed star},
  {magenta,mark=|}  }
\usepgfplotslibrary{fillbetween}
\usepackage[T1]{fontenc}
\usepackage{mathptmx}
\usepackage[scaled=.92]{helvet}
\usepackage{courier}

\usepackage[a4paper,left=2.5cm,right=2.5cm,top=3cm,bottom=2.5cm]{geometry}
\usepackage{authblk}
\usepackage{amsmath,bm}
\usepackage{upgreek}
\usepackage{mathtools}
\usepackage{xcolor}
\usepackage{algorithm}
\usepackage{algpseudocode}
\theoremstyle{remark}

\title{Generalised Adaptive Cross Approximation for Convolution
  Quadrature based Boundary Element Formulation}

\author[1]{A.M. Haider}
\author[2]{S. Rjasanow}
\author[1]{M. Schanz}
\affil[1]{Institute of Applied Mechanics, Graz University of
  Technology, Technikerstraße 4/II, 8010 Graz, Austria,
  m.schanz@tugraz.at}
\affil[2]{Department of Mathematics, Saarland University, 66041
  Saarbrücken, Germany, rjasanow@num.uni-sb.de}
\date{}                     %% if you don't need date to appear
\definecolor{facecol}{rgb}{0.16, 0.4, 1.}
\definecolor{fibercol}{rgb}{0.,0.65, 0.31}
\DeclareMathOperator*{\argmax}{arg\,max}
\newcommand{\op}[1]{\mathcal{#1}}
\newcommand{\kl}[1]{\left(#1\right)}
\newcommand{\Dt}{\Delta t}
\newcommand{\Od}[1]{\operatorname{d}#1}
\newcommand{\I}{\mathrm{i}}
\newcommand{\vek}[1]{\mathbf{#1}}      %% vectors (latin)
\newcommand{\x}[0]{\vek{x}}            %% position x
\newcommand{\y}[0]{\vek{y}}            %% position y
\newcommand{\mat}[1]{\mathsf{#1}}      %% matrices
\newcommand{\cvek}{\mat{c}}
\newcommand{\bvek}{\mat{b}}
\newcommand{\bone}{\mathbbm{1}}
\newcommand{\A}{\mat{A}}
\newcommand{\bt}{\mat{b}^{\mathsf{T}}}
\newcommand{\matWeightO}[1]{\mat{W}^{\Dt_n}_{#1}}
\newcommand{\matWeight}[2]{\matWeightO{#1}\!\!\kl{#2}}
\newcommand{\ie}{i.e., }

\newcommand{\eg}{e.g., }

\newcommand{\corr}[1]{{\color{black}#1}}
\graphicspath{{pics/}{pics_tmp}}
\begin{document}
	
\maketitle
\section*{Abstract}
    The acoustic wave equation is solved in time domain with a boundary
  element formulation. The time discretisation is performed with the
  generalised convolution quadrature method and for the spatial approximation
  standard lowest order elements are used. Collocation and Galerkin
  methods are applied. In the interest
  of increasing the efficiency of the boundary element method,  a
  low-rank approximation such as the adaptive cross approximation
  (ACA) is carried  out. We discuss about a generalisation of the ACA
  to approximate a three-dimensional array of data, \ie usual boundary
  element matrices at several complex frequencies. This method
  is used within the generalised convolution quadrature (gCQ) method to
  obtain a real time domain formulation. The behaviour of the proposed
  method is studied with three examples, a unit cube, a unit cube with
  a reentrant corner, and a unit ball. The properties of the method
  are preserved in the data   sparse representation and a significant
  reduction in storage is obtained. 

\textbf{Keywords}: wave equation; boundary element method; generalised
convolution quadrature; multivariate adaptive cross approximation

\section{Introduction}
% Introduction
Wave propagation problems appear often in engineering, \eg for non-destructive testing or exploring the underground. Any of such problems are formulated with hyperbolic partial differential equations, \eg in acoustics or elastodynamics. Sometimes the physical model allows to have a description with parabolic partial differential equations, \eg for thermal problems. Despite that mostly a linear theory is sufficient, the handling of space and time requires expensive discretisation methods and for scattering problems even an unbounded domain has to be considered.

There are several numerical methods for solving such space and time dependent problems. Methods exist based on finite differences for both variables, as well as finite element methods in space in combination with finite differences in time or discontinuous Galerkin methods. These domain-based methods are often very efficient but have difficulties to handle unbounded domains, \eg for scattering problems. The main problem is to truncate the computational domain with some approximations of the radiation condition. Especially for the latter class of problems, time domain boundary element methods (BEM) are a well established alternative because they inherently fulfil the radiation condition. The basis are
boundary integral equations by the use of retarded potentials as counterpart to the governing hyperbolic
partial differential equation. The mathematical theory goes back to
the beginning of the last century by Fredholm for scalar problems like acoustics and later by Kupradze~\cite{kupradze79} for vectorial problems in elasticity.  The mathematical background of time-dependent boundary integral equations is summarised by Costabel~\cite{costabel04} and extensively discussed in the textbook by Sayas~\cite{sayas2016}.

The first numerical realisation of a time domain boundary element formulation is originated by  Mansur \cite{man} in the 80th of the last century. Despite often used, this approach suffers from instabilities (see, \eg\cite{peirce97}). A stable space-time formulation has been published by Bamberger and Ha-Duong~\cite{bambduong}, which has been further explored by the group of Aimi~\cite{aimi08a,aimi12a}. These approaches work directly in time domain, whereas a transformation to Laplace- or Fourier-domain results in suitable formulations~\cite{cruri}  as well. These formulations in transformed domains have the advantage of being efficient from a memory point of view, since only elliptic problems have to be solved, which are of the size of one time step of the time-domain methods. However, to obtain the time-dependent solution several frequency-dependent solutions have to be computed, which are difficult to be solved with iterative solvers, especially for higher frequencies. Furthermore, suitable and mostly not physically motivated parameter choices have to be found for the inverse transformation technique. In contrast, the time-domain methods require a huge amount of memory since the matrices have to be stored for each time step (up to a cut-off in 3D). Somehow in between transformation and time-domain methods are BE formulations based on the convolution quadrature (CQ) method proposed by
Lubich~\cite{lubich88a,lubich88b}. Such a formulation is a true time stepping method utilising the Laplace domain fundamental solutions and properties. Applications of the CQ to BEM can be found, \eg in~\cite{schanz01a,schanz97e}. The
generalisation of this seminal technique to variable time step sizes
has been proposed by L\'opez-Fern\'andez and
Sauter~\cite{lopezfernandez13a, lopezfernandez15b} and is called
generalised convolution quadrature method (gCQ). Applications can be found in acoustics with absorbing boundary conditions~\cite{sauterschanz17a} and in thermoelasticity~\cite{leitner20a}. A comparison of a transformation and a CQ based BEM can be found in~\cite{schanz15a}.

The drawback of all BE formulations either for elliptic and much stronger for hyperbolic problems is the high storage and computing time demand as a standard formulation scales with $\mathcal{O}(M^2)$ for $M$ unknowns. In time domain, additionally, the time complexity has to be considered, where in the case of a CQ based formulation the complexity is of order $\mathcal{O}(M^2 N)$ for $N$ time steps. For elliptic problems fast methods has been proposed as the fast multipole method (FMM)~\cite{greengard97a} or $\mathcal{H}$-matrix based methods with the adaptive cross approximation (ACA) used in the matrix blocks~\cite{bebendorf03,bebendorf08a}. For hyperbolic problems the extension of FMM to the time variable has been published in~\cite{ergin98a} for acoustics and in~\cite{takahashi,otani06a} for elastodynamics.
In combination with CQ, fast methods are published in combination with a reformulation of CQ. The CQ-algorithm can be rewritten such that several elliptic problems have to be solved~\cite{banjai09a,schanz10a}. Based on this transformation the known techniques from the elliptic problems can be applied~\cite{banjai14a,messner10a}. 

Here,  a different approach is used. Independently whether the CQ in the original form or gCQ is used, essentially, a three-dimensional data array has to be efficiently  computed and stored. This data array is constructed by the spatial discretisation, resulting in two-dimensional data, and the used complex frequencies of the algorithm, which gives the third dimension. To find a low rank representation of this three-dimensional tensor the generalised adaptive cross approximation (3D-ACA) can be used. This technique is a generalisation 
of ACA~\cite{bebendorf03} and is proposed by
Bebendorf et al.~\cite{bebndorf11a,bebendorf13a}. It is based on a Tucker
decomposition~\cite{tucker66a} and can be traced back to  the group
of Tyrtyshnikov~\cite{oseledets08a}. The 3D-ACA might be seen as a higher-order singular value decomposition or as a multilinear SVD~\cite{lathauwer00a}, which is a generalisation of the matrix SVD to tensors.
An alternative but related approach is presented in~\cite{dirckx22a}. This approach does not look for a low rank representation but uses an interpolation with respect to the frequencies to reduce the total amount of necessary frequencies to establish the time domain solution.

Here, the 3D-ACA is applied on a gCQ based time domain formulation utilising the original idea of the multivariate ACA~\cite{bebendorf13a}. Recently a very similar approach has been published by Seibel~\cite{seibel22a}, where the conventional convolution quadrature method is used and, contrary to the approach proposed here, the $\mathcal{H}^2$ technique is used for the approximation of the boundary element matrices at the different frequencies.
\section{Problem statement}
% Problem statement
Let $\Omega \subset \mathbb{R}^3$ be a bounded Lipschitz domain and
$\Gamma = \partial\Omega$ its boundary with the outward normal
$\mathbf{n}$. The acoustic wave propagation is governed by
\begin{align} \label{eq:gov}
  \frac{\partial^2}{\partial t^2} u(\x,t)-
  c^2\Delta u(\x,t)& =0
                             \quad &(\x,t)&\in\Omega\times (0,T)\\
  u(\x,0)
  =\frac{\partial}{\partial t} u(\x,0)
                           &= 0
                             \quad & \x &\in \Omega \\
  \gamma_0 u(\x,t) &= g_D (\x,t) 
                    \quad & (\x,t) &\in \Gamma_D\times(0,T) \\
  \gamma_1 u(\x,t) &= g_N (\x,t) 
                    \quad & (\x,t) &\in \Gamma_N\times(0,T) 
\end{align}
with the wave speed $c$ and the end time $T>0$. $\Gamma_D$ and
$\Gamma_N$ denote either a Dirichlet or Neuman
boundary. The respective traces are defined by
\begin{align} \label{eq:traceD}
  \gamma_0 u(\x,t) = \lim_{\Omega \owns \x \to
  \x \in \Gamma} u(\x,t)
\end{align}
for the Dirichlet trace and with
\begin{align} \label{eq:traceN}
  \gamma_{1,x} u(\x,t) = \lim_{\Omega \owns \x \to
  \x \in \Gamma} \kl{\nabla u(\x,t) \cdot \mathbf{n}} = q(\x,t)
\end{align}
for the Neumann trace, also known as the conormal derivative or flux
$q(\x,t)$. 

The given problem can be solved with integral equations. Essentially,
either an indirect approach using layer potentials or an direct
approach are suitable. In the subsequent numerical tests both
metho\-dologies are used to show distinct features of the proposed
method. First, the retarded boundary integral operators have to be
formulated. It is the single layer potential 
\begin{align} \label{eq:slp}
  (\op{V} * q)(\x,t) = 
  \int_0^t \int_\Gamma \gamma_0U(\x - \y, t - \tau) q(\y, \tau) \Od{s}_\y
  \Od{\tau} \; ,
\end{align}
where the fundamental solution is $U(\x - \y, t - \tau) = \frac{1}{4 \pi \|\x-\y\|} \delta\kl{t-\tau -
  \frac{\|\x-\y\|}{c}}$ introducing a weak singularity. \corr{Here, $\delta(\cdot)$
denotes the Dirac distribution.} The double layer potential is defined by
\begin{align} \label{eq:dlp}
  (\op{K} * u)(\x,t) = \int_0^t  \int_{\Gamma} \gamma_{1,y} U(\x - \y, t -
  \tau) u( \y, \tau) \Od{s}_\y \Od{\tau} 
\end{align}
and its adjoint
\begin{align} \label{eq:adlp}
  (\op{K}' * u)(\x,t) = \int_0^t  \gamma_{1,x}  \int_{\Gamma} U(\x - \y, t -
  \tau) u( \y, \tau) \Od{s}_\y \Od{\tau} \; .
\end{align}
The remaining operator to be defined is the hypersingular operator
\begin{align} \label{eq:hso}
  (\op{D} * u)(\x,t) = - \int_0^t  \gamma_{1,x}  \int_{\Gamma}\hspace{-2.7ex} = \; \gamma_{1,y} U(\x - \y, t -
  \tau) u( \y, \tau) \Od{s}_\y \Od{\tau} \; ,
\end{align}
where $\int\hspace{-1.7ex} =$ indicates that this integral  has to
  be defined in the sense of Hadamard.
Using these operators, which are also called retarded potentials, the
solution for the Dirichlet problem related to \eqref{eq:gov} can be
obtained by 
\begin{equation} \label{eq:bie_diri}
  \kl{\op{V} \ast q} \!\!\kl{\x,t}  = \Bigl(\op{C}g_D +
    \kl{\op{K} \ast g_D}\Bigr)\!\!\kl{\x,t}
  \quad        (\x,t)\in\Gamma_D\times (0,T) 
\end{equation}
with $\Gamma = \Gamma_D$. On the right hand side the so-called
integral free term
\begin{align} \label{eq:cterm}
  \op{C}g_D(\x,t) =  \lim_{\varepsilon \to 0}   \int_{\partial
  B_\varepsilon(\x) \cap \Omega} \nabla \frac{1}{\|\x - \y\|} \cdot
  \mathbf{n} \; g_D( \y, t)
                               \Od{s}_\y + g_D( \x, t) \; .
\end{align}
has to be considered with $B_\varepsilon(\x)$ denoting  a ball of
radius $\varepsilon$ centered at $\x$ and $\partial
B_\varepsilon(\x)$ is its surface. In case of the Galerkin formulation
this expression reduces to $1/2 g_D(\x,t)$ but in the
collocation schema it is dependent on the geometry at the collocation
node. \corr{See, \eg~\cite{mantic93} for a description how the integral free
term depends on the geometry.} Alternatively to \eqref{eq:bie_diri}, the indirect approach 
\begin{equation} \label{eq:bie_indirect}
  (\op{V} * \vartheta)(\x,t) = g_D\!\kl{\x,t}  \quad (\x,t)\in\Gamma_D\times (0,T)
\end{equation}
can be used, where $\vartheta (\x,t)$ denotes the density
function. The Neumann problem can be solved with the second boundary integral equation
\begin{equation} \label{eq:bie_neum}
  \kl{\op{D} \ast u} \!\!\kl{\x,t}  = - \kl{\frac{1}{2}\op{I}g_N -
    \kl{\op{K}' \ast g_N}}\!\!\kl{\x,t}   \quad        (\x,t)\in\Gamma_N\times (0,T) \;,
\end{equation}
where the integral free term is already substituted by the expression
holding only for the Galerkin approach. The hypersingular
  integral is reduced to a weak singular one by partial integration,
  see, \eg~\cite{steinbach08a}.
\section{Boundary element method: Discretisation}
% Discrete integral equations

\paragraph{Spatial discretisation}
The boundary $\Gamma$ is discretised with elements resulting in an
approximation $\Gamma^h$ which is the union of geometrical elements
\begin{equation}
  \label{eq:geomapprox}
  \Gamma^h = \bigcup_{e=1}^{E} \tau_e\,.
\end{equation}
$\tau_e$ denote boundary elements, here linear surface triangles, and
their total number is $E$. 

Next, finite element bases on boundaries $\Gamma_D$ and $\Gamma_N$ are
used to construct the approximation spaces
\begin{align} \label{eq:spaces}
    X_D
    &= \operatorname{Span}\{\varphi_1,\varphi_2,\dots,\varphi_{M_1} \},\\
    X_N
    & =  \operatorname{Span}\{\psi_1,\psi_2,\dots,\psi_{M_2} \}.
\end{align}
The unknowns $u$ and  $q$ are approximated by a 
linear combination of functions in $X_D$ and $X_N$
\begin{equation}\label{eq:ansatz_space}
  u^h = \sum_{\ell=1}^{M_1} u_\ell(t) \varphi_\ell(\y) 
 \quad \text{and} \quad
 q^h = \sum_{k=1}^{M_2} q_k(t) \psi_k(\y)
 \quad \text{and} \quad
  \vartheta^h = \sum_{k=1}^{M_2} \vartheta_k(t) \psi_k(\y) \; .
\end{equation}
Note, the coefficients $u_\ell(t), q_k(t)$, and $\vartheta_k(t)$ are still continuous functions of
time $t$. In the following, the shape functions
$\varphi_\ell$ will be chosen linear continuous and $\psi_k$ constant
discontinuous. This is in accordance with the necessary function
spaces for the integral equations, which are not discussed here
(see~\cite{sayas2016}).

Using the spatial shape functions in the boundary integral
equations of the previous section results in semi-discrete
equations. They can be solved either with a Galerkin or a
collocation approach. Here, both techniques are presented to show
their respective behavior. In particular, the Neumann problem is solved
with a Galerkin method since in this case it is superior with respect
to the convergence behavior compared to the collocation method.

Hence, inserting the spatial shape functions from above in  the boundary integral
equations \eqref{eq:bie_diri} and applying the
collocation  method \corr{with the collocation points $\x_i$ at the nodes} results in the semi-discrete equation system
\begin{equation}  \label{eq:bie_semi_diri}
  \begin{split}
    \mat{V} \ast \mat{q}^h & = \mat{C} \mat{g}_D + \mat{K} \ast
    \mat{g}_D  \quad \text{with} \\
    \mat{V}[i,j] & = \int\limits_{\operatorname{supp} \psi_j}  \gamma_0 U(\corr{\x_i} - \y, t -
    \tau) \psi_j(\y)   \Od{s}_\y \qquad
    \mat{K}[i,j] = \int\limits_{\operatorname{supp} \varphi_j}  \gamma_{1,y} U(\corr{\x_i} - \y, t -
    \tau) \varphi_j(\y)   \Od{s}_\y 
  \end{split}
\end{equation}
for the Dirichlet problem. The alternative solution is obtained by
inserting the ansatz in \eqref{eq:bie_indirect} and applying the
Galerkin scheme
\begin{equation}  \label{eq:bie_semi_indirect}
  \mat{V} \ast \upvartheta^h = \mat{g}_D \qquad \text{with} \quad
  \mat{V}[i,j] = \int\limits_{\operatorname{supp} \psi_i} \psi_i(\x)
  \int\limits_{\operatorname{supp} \psi_j}  \gamma_0 U(\x - \y, t -
  \tau) \psi_j(\y)   \Od{s}_\y \Od{s}_\x
\end{equation}
The Neumann problem is solved by \eqref{eq:bie_neum} using these shape
functions and the Galerkin method resulting in 
\begin{equation}  \label{eq:bie_semi_neum}
  \begin{split}
    \mat{D} \ast \mat{u}^h  &= - \frac{1}{2} \mat{g}_N + \mat{K}' \ast
    \mat{g}_N
    \quad \text{with} \\
    \mat{D}[i,j]  & = \!\! \int\limits_{\operatorname{supp} \varphi_i}
    \varphi_i(\x) \gamma_{1,x} \!\! \int\limits_{\operatorname{supp} \varphi_j}
    \hspace{-4ex} = \; \gamma_{1,y} U(\x - \y, t -
    \tau) \varphi_j(\y)   \Od{s}_\y \Od{s}_\x\\ %\quad
    \mat{K}'[i,j] &= \!\! \int\limits_{\operatorname{supp} \varphi_i}
    \varphi_i(\x) \gamma_{1,x}\!\! \!\!\int\limits_{\operatorname{supp} \varphi_j}  \!U(\x - \y, t -
    \tau) \varphi_j(\y)   \Od{s}_\y  \Od{s}_\x .
  \end{split}
\end{equation}
The given boundary data are approximated by the same shape functions,
respectively. The notation with sans serif 
fonts or upright greek fonts in \eqref{eq:bie_semi_diri}, \eqref{eq:bie_semi_indirect}, and
\eqref{eq:bie_semi_neum} indicates that in these vectors the nodal
values are collected.

In the Galerkin formulations double integrations have to be performed,
whereas in the collocation schema only one spatial integration is
necessary. The weak singular integrals in the Galerkin
  formulations are treated with the formula by Erichson and
  Sauter~\cite{erichsen98} and in the collocation scheme with a Duffy
  transformation~\cite{duffy82a}. The regular integrals are treated
  with a standard Gaussian quadrature using a heuristic distance-based formula 
  to determine the number of Gauss points. No further treatment of
  quasi-singular integrals is considered.

\paragraph{Temporal discretisation}
The above given semi-discrete integral equations are discretised
in time using gCQ~\cite{lopezfernandez13a, lopezfernandez15b}. Here, the
variant using Runge-Kutta methods as the underlying time stepping
technique is applied. Details and the analysis can be found
in~\cite{lopezfernandez15b}, where here only the idea is briefly
sketched using the integral equation \eqref{eq:bie_semi_indirect}.

To present the idea of the gCQ, the convolution
  integral in \eqref{eq:bie_semi_indirect} is used. However, only
  one element of the spatial matrix vector product is shown to
  simplify notation. This convolution integral is
reformulated using the definition of the inverse Laplace transform
\begin{align} \label{eq:conv}
  (\mat{V}[i,j] \ast \upvartheta^h[j])(t) = \int_0^t 
  \mat{V}[i,j](t - \tau) \upvartheta^h[j](\tau) \Od{\tau} =
  \frac{1}{2 \pi \I} \int_C \hat{\mat{V}}[i,j](s) \int_0^t e^{s (t-\tau)}
  \upvartheta^h[j](\tau) \Od{\tau} \Od{s} \, ,
\end{align}
where for the Laplace variable holds $s \in \mathbb{C}, s.t. \Re{s} >0$. The symbol
$\hat{()}$ denotes the Laplace transform of a function, \corr{\eg
$\hat{\mat{V}}[i,j](s)$ means that in \eqref{eq:bie_semi_indirect} the
fundamental solution $U(\x - \y, t -
  \tau)$ is replaced by $\hat{U}(\x - \y, s)$}. Furthermore, the
imaginary unit is denoted by $\I=\sqrt{-1}$. The integration contour
$C$ is a vertical line $a-\I \infty \to a+\I \infty$, where $a$ is
larger than the real part of all singularities of the integral kernel. The rearrangement in
\eqref{eq:conv} is valid only if the Laplace transform of the kernel function and its inverse exist. This holds true for
the fundamental solutions above, where the Laplace transform of $U(\x-\y, t)$ is
$\hat{U}(\x-\y,s) = \frac{1}{4 \pi \|\x-\y\|} e^{-s\frac{\|\x-\y\|}{c}}$.  The inner time
integral is the solution of the differential equation of first order
\begin{align}\label{eq:dgl}
  \frac{\partial}{\partial t} x(t,s) = s x(t,s) + \upvartheta^h[j](t) \quad \text{with} \;\;
  x(t=0,s) =0 \;.
\end{align}
Thus, vanishing initial conditions are required. This ordinary differential equation can
be solved numerically by a time stepping method.  
Time is discretised in $N$ not necessarily constant time steps
$\Dt_i$, \ie 
\[ [0,T]=[0, t_1, t_2, \ldots, t_N], \quad  \Dt_i=t_i - t_{i-1}, \;\; 
  i=1,2, \ldots, N \, . \] 
Suppose an A- and L-stable Runge-Kutta
method is given by its Butcher tableau
\begin{tabular}{c|c}
  $\cvek$ & $\A$\\\hline    & $\bt$ 
\end{tabular}
with $\A \in \mathbb{R}^{m \times m}$, $\bvek,\cvek \in \mathbb{R}^m$
and $m$ is the number of stages. With the vector
$\bone=(1,1,\ldots,1)^{\mathrm{T}}$ of size $m$  the stability function is 
\begin{align}
  R(z) := 1 + z \bt (\mat{I} - z \A)^{-1} \bone 
\end{align}
and \corr{A- and L-stability} requires
\begin{align}
  |R(z)| \le 1 \quad \text{and} \quad
  (\mat{I} - z \A) \quad \text{is non-singular for} \; \Re{z} \le 0
  \quad \text{and} \quad
  \corr{\lim_{\Re{z} \to -\infty} R(z) = 0}
\end{align}
\corr{These conditions required that $\bt \A^{-1} =
(0,0,\ldots,1)$ holds and they imply $c_m=1$.} These severe restrictions on the used time
stepping method are  based partly on the proofs, but mostly on experience~\cite{schanz01a}.
With these definitions the discrete solution   of the ODE~\eqref{eq:dgl} is~\cite{lopezfernandez15b}
\begin{align}
  \mat{x}_n(s) = \kl{\mat{I} - \Dt_n s \A}^{-1} \kl{\kl{\bt \A^{-1} \mat{x}_{n-1}(s)} \bone +
  \Dt_n \A (\upvartheta^h[j])_n} \; .
\end{align}
Note that the solution $\mat{x}_n(s) $ is a vector containing $m$ results at the stages for the time step
$n$. The notation $()_n$ indicates the discrete \corr{values} of the
respective function/vector or of a product at \corr{all stages corresponding to time
$t_n$. Replacing in \eqref{eq:conv} the inner integral on
  the right hand side by its discrete approximation at all stages $m$ gives}
\begin{align}\label{eq:conv_discrete}
  \begin{split}
    \corr{(\mat{V}[i,j] \ast \upvartheta^h[j])_n} &= \frac{1}{2 \pi \I} \int_C \hat{\mat{V}}[i,j](s)
    \kl{\mat{I} - \Dt_n s \A}^{-1} \kl{\Dt_n \A (\upvartheta^h[j])_n +
      \kl{\bt \A^{-1} \mat{x}_{n-1}(s)} \bone           } \Od{s} \\
    & = \hat{\mat{V}}[i,j]\kl{\kl{\Dt_n \A}^{-1}} (\upvartheta^h[j])_n
    + \frac{1}{2 \pi \I} \int_C \hat{\mat{V}}[i,j](s) \kl{\mat{I} - \Dt_n s \A}^{-1}
    \kl{\bt \A^{-1} \mat{x}_{n-1}(s)} \bone \Od{s} \; .
  \end{split}
\end{align}
\corr{Note, the notation $(\mat{V}[i,j] \ast \upvartheta^h[j])_n$
  indicates a vector of size $m$.} The first part on the right-hand side is obtained by the residue theorem, while the
second term will be approximated by a quadrature rule using $N_Q$ integration
points and integration weights $\omega_\ell$.  The parameters and
integration weights used for the quadrature 
rule are given in~\ref{app:gcq}.  Approximating the convolution     in \eqref{eq:bie_semi_indirect} as above
gives the formula 
\begin{align} \label{eq:gcq_slp}
    \corr{(\mat{g}_D[i])_n} = \hat{\mat{V}}[i,j]\kl{\kl{\Dt_n \A}^{-1}} (\upvartheta^h[j])_n
    + \sum\limits_{\ell=1}^{N_Q} \omega_{\ell}
    \hat{\mat{V}}[i,j]\kl{s_{\ell}} \kl{\bt \A^{-1} \cdot \mat{x}_{n-1}(s_{\ell})} \kl{\mat{I} - \Dt_n
      s_{\ell} \A}^{-1} \bone \; .
\end{align}
The notation
$\hat{\mat{V}}[i,j]\kl{\kl{\Dt_n \A}^{-1}}$ means that the discrete operator is 
not evaluated at a single complex frequency as in $\hat{\mat{V}}[i,j]\kl{s_\ell}$ but at all
values of the matrix entries of $\kl{\Dt_n \A}^{-1}$. This
  matrix entry is therefore a matrix of size $m \times m$. Note,
$(\upvartheta^h[j])_n$ and $(\mat{g}_D[i])_n $ are also vectors of size $m$.

Eq \eqref{eq:gcq_slp} shows the computation of one entry in
  \eqref{eq:bie_semi_indirect} \corr{at all stages $m$}, which is the contribution of one
  node. Collecting these entries of all nodes in matrices
  $\hat{\mat{V}}(\cdot)$ and analogous defined vectors
  $\mat{g}_D(t_n)$ and $(\upvartheta^h)_n$, the algorithm can be given
  in a compact form. Note,
  the symbol $\hat{()}$ has still the meaning that the Laplace
  transform is applied on each entry and if the argument is the matrix
  $\kl{\Dt_n \A}^{-1}$ it means that in each entry we get a matrix of
  size $m \times m$. Considering this notation the
whole algorithm can be given as in~\cite{lopezfernandez12a} 
\begin{itemize}
\item First step
  \begin{equation*}
    \corr{(\mat{g}_D)_1 } = \hat{\mat{V}}\kl{\kl{\Dt_1 \A}^{-1}} (\upvartheta^h)_1
  \end{equation*}
  with implicit assumption of zero initial condition.
\item For all steps $n=2,\ldots,N$ the algorithm has two sub-steps
  \begin{enumerate}
  \item Update the solution vector $\mat{x}_{n-1}$ at all integration points $s_{\ell}$ 
    \begin{equation} \label{eq:gCQ_xn}
      \mat{x}_{n-1}\kl{s_{\ell}} = \kl{\mat{I} - \Dt_{n-1} s_{\ell} \A}^{-1}
      \kl{\kl{\bt \A^{-1} \cdot \mat{x}_{n-2}(s_{\ell})} \bone
        + \Dt_{n-1} \A (\upvartheta^h)_{n-1} }
    \end{equation}
    for $\ell=1,\ldots,N_Q$ with the number of integration points
    $N_Q$.
  \item Compute $(\upvartheta^h)_n $ solving the following linear system
    \begin{equation} \label{eq:gcq_formula}
        \corr{(\mat{g}_D)_n} =  \hat{\mat{V}}\kl{\kl{\Dt_n \A}^{-1}}
        (\upvartheta^h)_n + \sum\limits_{\ell=1}^{N_Q}
        \hat{\mat{V}}\kl{s_\ell}   \matWeight{\ell}{(\upvartheta^h)_{n-1}}
        \; .
      \end{equation}
      The abbreviation $\matWeight{\ell}{(\upvartheta^h)_{n-1}}$
      denotes a vector of size $M_2 m$, which is dependent on
      $\mat{x}_{n-1}(s_{\ell})$. It collects the entries $\omega_{\ell} \kl{\bt \A^{-1} \cdot
        \mat{x}_{n-1}(s_{\ell})} \kl{\mat{I} - \Dt_n s_{\ell} \A}^{-1}
      \bone$  in the sum in \eqref{eq:gcq_slp} related to each node.
  \end{enumerate}
\end{itemize}
Essentially, this algorithm requires the evaluation of the integral kernel at $N_Q$ points $s_\ell$,
which are complex frequencies. Consequently, we get an array of system matrices of
  size $M_2 \times M_2 \times N_Q$ beside the matrix $ \hat{\mat{V}}\kl{\kl{\Dt_n \A}^{-1}}$.
This array of system matrices can be interpreted as a three-dimensional
array of data which will be approximated by a data-sparse representation based on 3D-ACA.

Before describing this algorithm, the discrete set of integral
equations is given. For the Dirichlet problem, we have the direct approach
\begin{equation}  \label{eq:bie_discr_diri}
  \hat{\mat{V}}\kl{\kl{\Dt_n \A}^{-1}} (\mat{q}^h)_n = \kl{\mat{C} +
    \hat{\mat{K}}\kl{\kl{\Dt_n \A}^{-1}} } (\mat{g}_D)_n + \sum\limits_{\ell=1}^{N_Q} \kl{
    \hat{\mat{K}}\kl{s_\ell} \matWeight{\ell}{(\mat{g}_D)_{n-1}} - \hat{\mat{V}}\kl{s_\ell}
  \matWeight{\ell}{(\mat{q}^h)_{n-1}} }
\end{equation}
and the indirect approach
\begin{equation} \label{eq:bie_discr_indirect}
  \hat{\mat{V}}\kl{\kl{\Dt_n \A}^{-1}} (\upvartheta^h)_n  =
  (\mat{g}_D)_n -  \sum\limits_{\ell=1}^{N_Q} \hat{\mat{V}}\kl{s_\ell}
  \matWeight{\ell}{(\upvartheta^h)_{n-1}} \; .
\end{equation}
The Neuman problem is solved with the discrete equations
\begin{equation}  \label{eq:bie_discr_neum}
  \begin{split}
    \hat{\mat{D}}&\kl{\kl{\Dt_n \A}^{-1}} (\mat{u}^h)_n =\\  &-
  \kl{\frac{1}{2} \mat{I} - \hat{\mat{K}}'\kl{\kl{\Dt_n \A}^{-1}}}(\mat{g}_N)_n + \sum\limits_{\ell=1}^{N_Q} \kl{
    \hat{\mat{K}}'\kl{s_\ell} \matWeight{\ell}{(\mat{g}_N)_{n-1}} - \hat{\mat{D}}\kl{s_\ell}
    \matWeight{\ell}{(\mat{u}^h)_{n-1}} } \; .
  \end{split}
\end{equation}
By examining these equations the computational cost can be estimated. The computation
of the matrices is $\mathcal{O}((N_Q+1) M^2)$ and the evaluation of the time stepping
method is $\mathcal{O}(N_Q N)$ matrix-vector multiplications. Note, as given
in~\ref{app:gcq} for Runge-Kutta methods with stages $m>1$ a suitable choice is $N_Q=N
(\log(N))^2$.  The solution of the equation systems in
\eqref{eq:bie_discr_diri}, \eqref{eq:bie_discr_indirect}, or \eqref{eq:bie_discr_neum} contributes with
$\mathcal{O}(M^2 n_{iter})$. An iterative solver is necessary as long as the
  matrix of the actual time step is approximated with FMM or
  $\mathcal{H}$-technique. In these estimates the spatial
dimension is denoted by $M$ which is either $M_2$ or $M_1$, respectively. The value $n_{iter}$
is the number of iterations for the equation solver, which is
usually small. 
\section{Three-dimensional adaptive cross approximation}
% 3D-ACA
An approximation of a three-dimensional  array of data or a tensor of
third order $\mathcal{C}\in\mathbb{C}^{M\times M\times N_Q}$
in terms of a low-rank approximation has been proposed in~\cite{bebendorf13a}
and is referred to as a generalisation of adaptive cross approximation
or also called  3D-ACA. In this approach, the 3D array of data to
be approximated is generated by defining the outer product by
\begin{equation} \label{eq:outerProd}
  \mathcal{C} = \mathbf{H} \otimes \mathbf{f}
\end{equation}
with $ \mathbf{H}\in\mathbb{C}^{M\times M},\  \mathbf{f}\in
\mathbb{C}^{N_Q}$. The matrix $\mathbf{H}$ corresponds to the spatial
discretization of the different potentials used in the boundary
element formulations from above at a specific frequency $s_{\ell}$,
\eg the single layer potential in \eqref{eq:bie_discr_indirect} on the
right hand side. This matrix will be called face or slice and may be computed as a
dense matrix or with a suitable hierarchical partition of the mesh as
$\mathcal{H}$-matrix. Within this matrix different low-rank
approximations can be applied. Here, the ACA~\cite{bebendorf03} is
used. Connected to this choice an hierarchical structure of the
matrix has to be build with  balanced cluster
trees~\cite{bebendorf08a}. The vector $\mathbf{f}$, called fiber in
the following, collects selected elements of $\mathbf{H}$  at the set
of frequencies determined by the gCQ. The 
latter would amount in $N_Q$ entries and the same amount of
faces. The 3D-ACA on the one hand approximates, as above mentioned,
the faces with low rank matrices and, more importantly, the amount of
necessary frequencies is adaptively determined. Hence, a sum of outer
products $\mathcal{C}^\ell$ as given in \eqref{eq:outerProd} is established. The
summation is stopped if  $\mathcal{C}^\ell$ is comparable to
$\mathcal{C}$ up to a prescribed precision $\varepsilon$. This is
measured with a Frobenius-norm.  The basic concept of this approach is sketched
in Algorithm~\ref{algo:3d_aca}.
\begin{algorithm}
  \caption{\label{algo:3d_aca}Pseudo code of 3D-ACA (taken from~\cite{seibel22a})}
  \begin{algorithmic}
    \Function{3D ACA}{$ENTRY, \varepsilon$}       \Comment{$ENTRY$
      provides the integrated kernel values at collocation
      point $x_i$ and element $j$}
    \State $\mathcal{C}^{(0)} =0, k_1=0$ and $\ell=0$
    \While{$\|\mathbf{H}_{\ell}\|_F \|\mathbf{f}_{\ell}\|_2 >
      \varepsilon \|\mathcal{C}^{(\ell)}\|_F$}
    \State $\ell = \ell+1$
    \State $H_{\ell}[i,j] = ENTRY(i,j,k_{\ell}) - \mathcal{C}^{(\ell -
      1)}[i,j,k_{\ell}], \qquad i,j = 1, \ldots , M$
    \State $H_{\ell}[i_{\ell},j_{\ell}] = \max_{i,j} | H_{\ell}[i,j]|$
    \State $f_{\ell}[k] = H_{\ell}[i_{\ell},j_{\ell}]^{-1} \kl{ENTRY(i_{\ell},j_{\ell},k) - \mathcal{C}^{(\ell -
        1)}[i_{\ell},j_{\ell},k]}, \qquad k = 1, \ldots , N_Q$
    \State $\mathcal{C}^{(\ell)} = \mathcal{C}^{(\ell-1)} +
    \mathbf{H}_{\ell}  \otimes \mathbf{f}_{\ell}$
    \State $k_{\ell+1} = \arg \max_k |f_{\ell}[k]|$
    \EndWhile
    \State $r = \ell -1$   \Comment{Final rank, \ie necessary
      frequencies}
    
    \Return $\mathcal{C}^r = \sum\limits_{\ell=1}^r \mathbf{H}_{\ell}  \otimes \mathbf{f}_{\ell}$
    \EndFunction
  \end{algorithmic}
\end{algorithm}
Essentially, three-dimensional crosses are established. These crosses
consist of a face $\mathbf{H}_{\ell}$, which is the respective matrix at a distinct
frequency $s_{\ell}$ times a fiber $\mathbf{f}_{\ell}$, which contains
one matrix entry, the pivot element,  at all $N_Q$ frequencies. These
crosses are added up until a predefined error is obtained. The amount
of used frequencies is the so-called rank $r$ of this approximation of
the initial data cube. The mentioned pivot element is denoted in
Algorithm~\ref{algo:3d_aca} with $H_{\ell}[i_{\ell}, j_{\ell}]$. Like
in the 'normal'  ACA this pivot element can be selected
arbitrarily. However, selecting the maximum element ensures
convergence. For details see~\cite{bebendorf13a} or the application
using $\mathcal{H}^2$-matrices in the faces~\cite{seibel22a}.

The stopping criterion requires a norm evaluation of
$\mathcal{C}^{(\ell)}$. Assuming some monotonicity 
the norm can be computed recursively
\begin{equation} \label{eq:norm}
  \left\| \mathcal{C}^{(\ell)}\right\|_F^2 = \sum\limits_{d,d'}^\ell \kl{\sum\limits_{i,j} H_d[i,j]
    \widebar{H_{d'}[i,j]}} \kl{\sum\limits_k f_d[k] \widebar{f_{d'}[k]}}
\end{equation}
following~\cite{seibel22a}.

The multiplication of the three-dimensional data array with a vector is changed by
the proposed algorithm. Essentially, the algorithm separates the
frequency dependency such that $\mathbf{H}_\ell\kl{\hat{\mat{V}}}$ is
independent of the frequency and the fibers
$\mathbf{f}_\ell\kl{\hat{\mat{V}}}$ present this dependency. Let us use the multiplication on the right hand side
of \eqref{eq:bie_discr_indirect} as example. The multiplication is changed to
\corr{\begin{align} \label{eq:product}
\sum\limits_{\ell=1}^{N_Q} \hat{\mat{V}}\kl{s_\ell}
  \matWeightO{\ell} =
  \sum\limits_{\ell=1}^{N_Q}  \sum\limits_{k=1}^r \mathbf{H}_k \kl{\hat{\mat{V}}}
  \otimes f_k[\ell]\!\!\kl{\hat{\mat{V}}}
  \matWeightO{\ell}
  = \sum\limits_{k=1}^r \mathbf{H}_k \kl{\hat{\mat{V}}}
  \otimes  \sum\limits_{\ell=1}^{N_Q} f_k[\ell]\!\!\kl{\hat{\mat{V}}}
  \matWeightO{\ell} \; .
      \end{align}
      }
The complexity of the original
operation is $\mathcal{O}\kl{N_Q M^2}$ for $M$ spatial unknowns. The approximated
version has the complexity $\mathcal{O}\kl{r (M^2 + N_Q)}$. It
consists of the inner sum, \corr{which requires  $N_Q$ multiplications}, and the outer sum which is a matrix
times vector multiplication of size $M$. Hence, the leading term with $M^2$ has
only a factor of $r$ instead of $N_Q$ compared to the dense
computation. It can be expected that for larger problem sizes with a
significant reduction from $N_Q$ to $r$ complex frequencies this
discrete convolution is faster.

% \noindent\textbf{Remarks for the implementation:}
\paragraph{Remarks for the implementation:}
\begin{itemize}
\item There are two possibilities how to apply the 3D-ACA
  within the above given boundary element formulations. Either it can be applied on
  each matrix block, which requires that for the different frequencies
  the same hierarchical structure is used. This holds here as a pure
  geometrical admissibility criterium for the $\mathcal{H}$-matrix has
  been selected. The other choice would be to apply 3D-ACA on the
  $\mathcal{H}$-matrix as a whole. In the latter case the stopping
  criterion is determined by all matrix blocks, \ie \eqref{eq:norm} is
  evaluated for each block and summed up. Consequently, the matrix block
  with the largest rank determines the overall rank. Hence, this choice
  results in a less efficient formulation as shown
  in~\cite{schanz23a}. Nevertheless, if the matrix structure in the
  different $\mathbf{H}_\ell$ differ this choice is a solution.
\item In each face the ACA is used to compress these matrices beside
  the rank reduction with respect to the frequencies. Here, a
  recompression of the ACA is applied, which gives a representation of
  the admissible matrix blocks in a SVD like form $\mathcal{A}  \approx
  \mathcal{A}_r = \mathbf{U} \Sigma \mathbf{V}^H$. Starting from the
  low rank representation of a block $\mathcal{A}_r = \mathbf{U}
  \mathbf{V}$, a QR-decomposition of the low-rank matrix $\mathbf{U}$
  and $ \mathbf{V}$ is performed 
  \begin{align}
    \mathcal{A}  \approx \mathcal{A}_r  = \mathbf{U} \mathbf{V}^H=
    Q_U( R_U R_V^H ) Q_V^H 
    = Q_U\,\check{\mathbf{U}}\check{\Sigma}\check{\mathbf{V}}^H Q_V^H
    = \overline{\mathbf{U}}\check{\Sigma}\overline{\mathbf{V}}^H
  \end{align}
  and, secondly, the SVD is applied on the smaller inner matrix $ R_U
  R_V^H$. The sketched technique improves
  the compression rates.
\item Within the algorithm the largest element in a matrix block has
  to be found. For non-admissible blocks this operation requires to
  compare all block entries. For admissible blocks the elements are
  not available directly and, hence, has either be computed out of the low-rank
  representation or an estimate can be used. Due to the
  recompression of the ACA-blocks low-rank matrices
  $\overline{\mathbf{U}}$ and $\overline{\mathbf{V}}$ consist of  
  orthonormal columns and the estimation of the index $i$ of the
  pivot position can be found with
  \begin{align} \label{eq:index}
    \argmax \limits_{i,j} \lvert a_{ij}\rvert =
    \argmax \limits_{i,j} \sum_{k=1}^r \lvert u_{ik}
    \sigma_k\rvert     \underbrace{\lvert     v_{jk}\rvert}_{\leq 1}
    \leq \argmax\limits_{i}  \sum_{k=1}^r \lvert u_{ik}
    \sigma_k\rvert  \; .
  \end{align}
  The analogous way is used also for index $j$. This algorithm gives
  a pivot element which is approximately the largest element. However,
  in the computation of the fiber the original entries are
  used. Hence, the computation of the residuum is an
  approximation. For the single layer potential this approximation
  works fine but not for the double layer potential. In this case, the
  maximal element is determined by multiplying the low rank matrices,
  finding the largest element and deleting this dense
  block. Fortunately, such matrix operations can be efficiently
  implemented such that the computing time is not affected. Hence, in
  the presented results for all potentials this expanding to the
  dense matrix block is used.
\end{itemize}

\section{Numerical examples}
% numerische Ergebnisse
% Zuerst tests am Wuerfel und dann an einspringender Ecke
% am Schluss Kugel
The above proposed method to accelerate the gCQ based time domain
boundary element method is tested at three examples. First, a unit
cube is used, second a cube with a reentrant corner
(three-dimensional L-Shape), and third a unit ball. The first two test
geometries are used for Dirichlet and Neumann problems, whereas in
case of the unit ball only a Dirichlet problem is studied. For the latter
case a non-smooth
timely response of the solution has to be tackled. In all tests the
main goal is to show that the approximation of the 3D-ACA does not
spoil the results, \ie the newly introduced approximation error is
smaller than the error of the dense BE formulation. Main focus is on
the overall obtained reduction in storage, \ie the compression which is
measured as the relation between the dense storage and the storage
necessary for the proposed method.

To show that the 3D-ACA does not spoil the results the error in space
and time is measured with
\begin{align} \label{eq:error}
  \begin{split}
    L_{max}(u_h)  &= \max\limits_{1 \le n \le N} \|
    u\kl{\frac{t_n+t_{n+1}}{2}} -  u_h\kl{\frac{t_n+t_{n+1}}{2}}
    \|_{L_2} \\ & = \max\limits_{1 \le n \le N} \sqrt{\int\limits_\Gamma 
      \kl{u\kl{\mathbf{x},\frac{t_n+t_{n-1}}{2}} - 
        u_h\kl{\mathbf{x},\frac{t_n+t_{n-1}}{2}} }^2 \Od{\Gamma} } \;
    ,
  \end{split}
\end{align}
where $u$ and $u_h$ are placeholders for the presented quantity
measured at the boundary of the selected domain $\Gamma$.
The same error has been used in~\cite{sauterschanz17a}. It computes
  the usual $L_2$-error with respect to the spatial variable and
  takes the largest value of them in all time steps.  The convergence
rate is denoted by 
\begin{equation} \label{eq:eoc}
  \mathrm{eoc} = \log_2\kl{\frac{L_{max}^h}{L_{max}^{h+1}} }
  \; ,
\end{equation}
where the superscript $h$ and $h+1$ denotes two subsequent refinement
level. The convergence results 
presented below does not change significantly 
if an $L_2$-error in time would be used as long as the time response
is smooth.  For the results of the sphere in the third subsection
the drop in convergence is not that significant for an $L_2$-error in
time as such a measure smears the error in time over all time steps.

As mentioned above, linear triangles are used for the discretisation
of the geometry and either linear continuous or constant
discontinuous  shape functions for the unknown Cauchy data in the
direct method are applied. In the indirect method constant discontinuous  shape functions for the
density are selected. The underlying time stepping
method has been the 2-stage Radau IIA. The approximation within the
faces has been chosen to be $\varepsilon_{ACA}=10^{-4} \ldots 10^{-8}$ for
the different refinement levels of the spatial discretisation. In each refinement
level the mesh size is halved as well as the time step size. The ratio
of time steps to mesh size was kept constant with 0.7. The precision
of the method with respect to the frequencies, \ie the $\varepsilon$
in  Algorithm~\ref{algo:3d_aca} was selected as $\varepsilon =
100*\varepsilon_{ACA}$. The final equation is solved with BiCGstab
without any preconditioner, where the precision $\varepsilon_{ACA}$ is
set. All results are computed for a wave speed $c=\unitfrac[1]{m}{s}$.

\subsection{Example: Unit cube loaded by a smooth pulse}
\begin{figure}[htb]
  \subcaptionbox{Unit cube (level 1, $h=\unit[0.5]{m}$)}{
    \includegraphics[scale=.23]{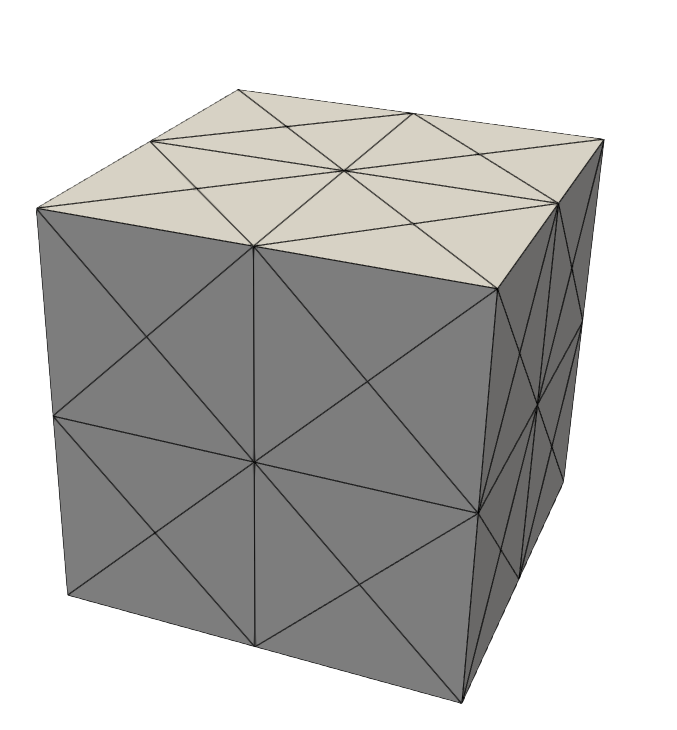}
  } \hfill
  \subcaptionbox{Used meshes}{
    \begin{tabular}{l|rrrr}
      level & nodes & elements & $h$ & $\Dt$\\ \hline
      1 & 50 & 96 & \unit[0.5]{m} & \unit[0.3]{s} \\
      2 & 194 & 384 & \unit[0.25]{m} & \unit[0.15]{s} \\
      3 & 770 & 1536 & \unit[0.125]{m} & \unit[0.075]{s} \\
      4 & 3074 & 6144 & \unit[0.0625]{m} & \unit[0.0375]{s} \\
      5 & 12290 & 24576 & \unit[0.03125]{m} & \unit[0.01875]{s}
    \end{tabular}
  }
  \caption{\label{fig:unitCube}Unit cube: Geometry and discretisation parameters}
\end{figure}

The first used test geometry is a unit cube $[-0.5, 0.5]^3$ with the
coordinate system located in the middle of the cube. In
Fig.~\ref{fig:unitCube}, this cube is displayed with the mesh of
refinement level 1 and the
table of all used meshes is given aside.
These meshes are created by bisecting the cathetus of the coarser
mesh. As load a smooth pulse
\begin{align} \label{eq:load}
  u\kl{\mathbf{y},t} = \frac{\kl{t-\frac{r}{c}}^2}{r}  e^{-c
  \kl{t-\frac{r}{c}}} \; H\kl{t-\frac{r}{c}} \qquad \text{with} \quad
  r=\|\mathbf{x}-\mathbf{y}\| 
\end{align}
at the excitation point $\mathbf{x}=(0.8, 0.2, 0.3)^{\mathsf{T}}$  is
used for the Dirichlet problem.  For the Neumann problem the same
geometry but with inverted normals, \ie a cavity with the shape of the
cube, is used. In this case the excitation point is at
$\mathbf{x}=(0.3, -0.2, 0.4)^{\mathsf{T}}$. The total observation time
$T=\unit[3]{s}$ is selected such that the smooth pulse travels over
the whole unit cube.

\begin{figure}[hbt]
  \centering
  % \begin{tikzpicture}
  %   \begin{loglogaxis}
  %     [scale=.9,
  %     cycle list name=mylines,
  %     scale only axis,
  %     xlabel=mesh $h$,
  %     xtick={0.03,0.06,0.12,0.25,0.5},
  %     ylabel= error $L_{max}$,
  %     log ticks with fixed point,
  %     legend columns=2,
  %     legend style={at={(.6,.2)}}
  %     ]
  %     \addlegendimage{empty legend};
  %     \addlegendentry{Dirichlet};
  %     \addlegendimage{empty legend};\addlegendentry{Neumann}; 
  %     %
  %     \addplot[red,dashed] table [x=h, y=Lmax_nolwr, col sep=comma] {data_cube/ErrorsDiriProblem.csv};
  %     \addlegendentry{No Lwr};
  %     \addplot[blue,dashed] table [x=h, y=Lmax_nolwr, col sep=comma] {data_cube/ErrorsNeumProblem.csv};
  %     \addlegendentry{No Lwr};
  %     \addplot[red,dotted]  table [x=h, y=Lmax_recaca, col sep=comma] {data_cube/ErrorsDiriProblem.csv};
  %     \addlegendentry{ACA};
  %     \addplot[blue,dotted] table [x=h, y=Lmax_recaca, col sep=comma] {data_cube/ErrorsNeumProblem.csv};
  %     \addlegendentry{ACA};
  %     \addplot[red, only marks] table [x=h, y=Lmax_dense, col sep=comma] {data_cube/ErrorsDiriProblem.csv};
  %     \addlegendentry{Dense};
  %     \addplot[blue, only marks] table [x=h, y=Lmax_dense, col sep=comma] {data_cube/ErrorsNeumProblem.csv};
  %     \addlegendentry{Dense};
  %     %
  %     \addplot[dashdotted] gnuplot[domain=.05:.3] {2.5*x};
  %     \addlegendentry{eoc = 1}; 
  %     \addplot[dashed] gnuplot[domain=.1:.35] {1.2*x*x};
  %     \addlegendentry{eoc = 2};
  %   \end{loglogaxis}
  % \end{tikzpicture}
  \includegraphics{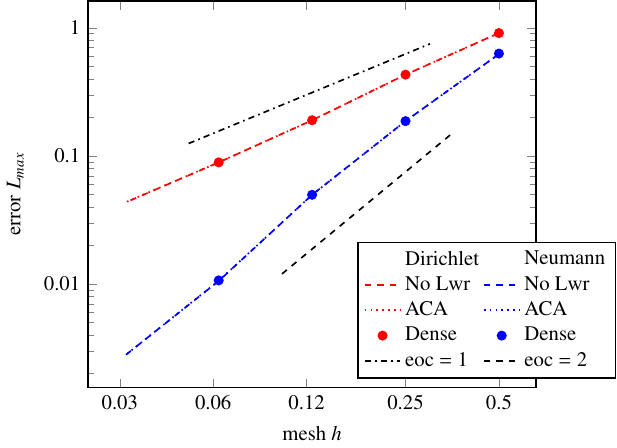}
  \caption{Cube: $L_{max}$-error versus refinement in space and time}
  \label{fig:ConvCube}
\end{figure}

In the following, the results are discussed for a Dirichlet problem using \eqref{eq:bie_discr_diri}
with a collocation approach and a Neumann problem using
\eqref{eq:bie_discr_neum} with a Galerkin approach. In
Fig.~\ref{fig:ConvCube}, the $L_{max}$-error defined in
\eqref{eq:error} is presented versus the different refinement levels,
where `No Lwr' means a computation applying the 3D-ACA but for the
faces instead of the ACA dense matrix blocks are applied, \ie the
approximation is only with respect to the complex frequencies.
The phrase `dense' means no 3D-ACA at all and `ACA' states that the
algorithm from above is used with an ACA approximation also in the faces.
Note, the refinement is done in space and time, \ie not only $h$ is
changed but as well $\Dt$.  The numbers displayed at
this axis refer to the element size $h$. The expected oder of
convergence is obtained considering that two factors influence the
error. First, the spatial discretisation is expected to yield for a
refinement in space a linear order for the constant shape functions
used for the single layer potential, whereas the linear shape
functions would result in a quadratic order for the hypersingular
operator. These spatial known erros for the respective elliptic
problem are combined with the 2-stage Radau IIA method used within the
gCQ. This method would allow a third order convergence. However, as
the spatial convergence orders are smaler they will dominate. This is
clearly observed in Fig.~\ref{fig:ConvCube}. The second observation
is that the 3D-ACA does not spoil the results. Certainly, the
$\varepsilon$-parameters have to be adjusted carefully. As mentioned
above, the ACA in the faces starts with $\varepsilon_{ACA}=10^{-4}$ in
the coarsest level and is increased by $10^{-1}$ in each level up to
$\varepsilon_{ACA}=10^{-8}$. The gCQ Algorithm~\ref{algo:3d_aca} is
terminated with $\varepsilon = 100*\varepsilon_{ACA}$.  Both
parameters are found by trial and error.

\begin{figure}[hbt]
  \subcaptionbox{Single layer potential}{
    % \begin{tikzpicture}
    %   \begin{loglogaxis}
    %     [scale=.6,
    %     scale only axis,
    %     xlabel=mesh $h$,
    %     ylabel= Number frequencies,
    %     log ticks with fixed point,
    %     legend columns=2,
    %     legend style={at={(0,1.1)}}
    %     ]
    %     \addplot [black,mark=square] table [x=h, y=no_f, col sep=comma] {data_cube/FreqsDiriProblem_nolwr.csv}; 
    %     \addlegendentry{No 3D-ACA};
    %     % 
    %     \addplot [name path = minfreq,blue] table [x=h, y=fmat_min, col sep=comma] {data_cube/FreqsDiriProblem_recaca.csv};
    %     \addlegendentry{Min no Freqs};
    %     \addplot [name path = maxfreq, blue] table [x=h, y=fmat_max, col sep=comma] {data_cube/FreqsDiriProblem_recaca.csv};
    %     \addlegendentry{Max no Freqs};
    %     \addplot [blue,mark=+] table [x=h, y=fmat_mean, col sep=comma] {data_cube/FreqsDiriProblem_recaca.csv};
    %     \addlegendentry{Mean no Freqs};
    %     \addplot [blue!20] fill between[of=maxfreq and minfreq]; 
    %     % 
    %     \addplot [red] table [x=h, y=fmat_min, col sep=comma] {data_cube/FreqsDiriProblem_nolwr.csv}; 
    %     \addplot [red] table [x=h, y=fmat_max, col sep=comma] {data_cube/FreqsDiriProblem_nolwr.csv};
    %     \addplot [red,mark=+] table [x=h, y=fmat_mean, col sep=comma]
    %     {data_cube/FreqsDiriProblem_nolwr.csv};
    %   \end{loglogaxis}
    % \end{tikzpicture}
    \includegraphics{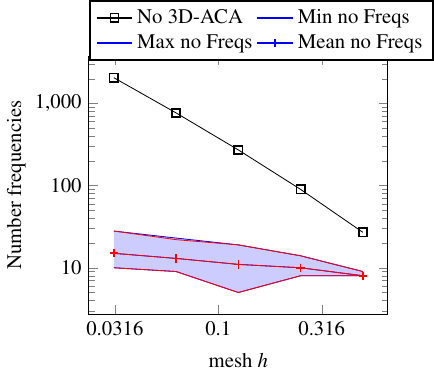}
  } \hfill
  \subcaptionbox{Double layer potential}{
    % \begin{tikzpicture}
    %   \begin{loglogaxis}
    %     [scale=.6, 
    %     scale only axis,
    %     xlabel=mesh $h$,
    %     log ticks with fixed point,
    %     legend style={at={(-.2,1.1)}}
    %     ]
    %     \addlegendimage{red} \addlegendentry{No lwr};
    %     \addlegendimage{blue} \addlegendentry{faces ACA};
    %     \addplot [black,mark=square] table [x=h, y=no_f, col sep=comma] {data_cube/FreqsDiriProblem_nolwr.csv}; 
    %     % 
    %     \addplot [name path = minfreq,blue] table [x=h, y=frhs_min, col sep=comma] {data_cube/FreqsDiriProblem_recaca.csv};
    %     \addplot [name path = maxfreq, blue] table [x=h, y=frhs_max, col sep=comma] {data_cube/FreqsDiriProblem_recaca.csv};
    %     \addplot [blue,mark=+] table [x=h, y=f_rhs_mean, col sep=comma] {data_cube/FreqsDiriProblem_recaca.csv};
    %     \addplot [blue!20] fill between[of=maxfreq and minfreq]; 
    %     % 
    %     \addplot [red] table [x=h, y=frhs_min, col sep=comma] {data_cube/FreqsDiriProblem_nolwr.csv}; 
    %     \addplot [red] table [x=h, y=frhs_max, col sep=comma] {data_cube/FreqsDiriProblem_nolwr.csv};
    %     \addplot [red,mark=+] table [x=h, y=f_rhs_mean, col sep=comma]
    %     {data_cube/FreqsDiriProblem_nolwr.csv};
    %   \end{loglogaxis}
    % \end{tikzpicture}
    \includegraphics{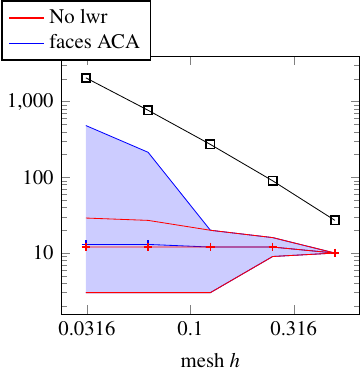}
  }
  \caption{\label{fig:CubeFreqsDiri}Cube: Number of used frequencies for the
    Dirichlet problem}
\end{figure}

\begin{figure}[hbt]
  \subcaptionbox{Hyper singular operator}{
    % \begin{tikzpicture}
    %   \begin{loglogaxis}
    %     [scale=.6,
    %     scale only axis,
    %     xlabel=mesh $h$,
    %     ylabel= Number frequencies,
    %     log ticks with fixed point,
    %     legend columns=2,
    %     legend style={at={(0,1.1)}}
    %     ]
    %     \addplot [black,mark=square] table [x=h, y=no_f, col sep=comma] {data_cube/FreqsNeumProblem_nolwr.csv}; 
    %     \addlegendentry{No 3D-ACA};
    %     % 
    %     \addplot [name path = minfreq,blue] table [x=h, y=fmat_min, col sep=comma] {data_cube/FreqsNeumProblem_recaca.csv};
    %     \addlegendentry{Min no Freqs};
    %     \addplot [name path = maxfreq, blue] table [x=h, y=fmat_max, col sep=comma] {data_cube/FreqsNeumProblem_recaca.csv};
    %     \addlegendentry{Max no Freqs};
    %     \addplot [blue,mark=+] table [x=h, y=fmat_mean, col sep=comma] {data_cube/FreqsNeumProblem_recaca.csv};
    %     \addlegendentry{Mean no Freqs};
    %     \addplot [blue!20] fill between[of=maxfreq and minfreq]; 
    %     % 
    %     \addplot [red] table [x=h, y=fmat_min, col sep=comma] {data_cube/FreqsNeumProblem_nolwr.csv}; 
    %     \addplot [red] table [x=h, y=fmat_max, col sep=comma] {data_cube/FreqsNeumProblem_nolwr.csv};
    %     \addplot [red,mark=+] table [x=h, y=fmat_mean, col sep=comma]
    %     {data_cube/FreqsNeumProblem_nolwr.csv};
    %   \end{loglogaxis}
    % \end{tikzpicture}
    \includegraphics{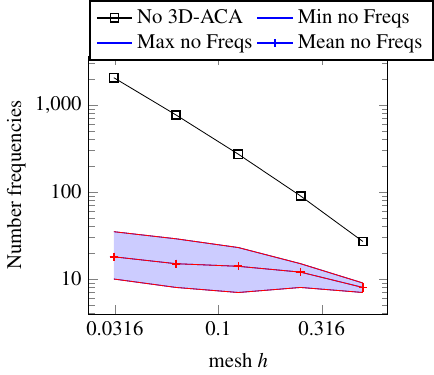}
  } \hfill
  \subcaptionbox{Adjoint double layer potential}{
    % \begin{tikzpicture}
    %   \begin{loglogaxis}
    %     [scale=.6, 
    %     scale only axis,
    %     xlabel=mesh $h$,
    %     log ticks with fixed point,
    %     legend style={at={(-.2,1.1)}}
    %     ]
    %     \addlegendimage{red} \addlegendentry{No lwr};
    %     \addlegendimage{blue} \addlegendentry{faces ACA};
    %     \addplot [black,mark=square] table [x=h, y=no_f, col sep=comma] {data_cube/FreqsNeumProblem_nolwr.csv}; 
    %     % 
    %     \addplot [name path = minfreq,blue] table [x=h, y=frhs_min, col sep=comma] {data_cube/FreqsNeumProblem_recaca.csv};
    %     \addplot [name path = maxfreq, blue] table [x=h, y=frhs_max, col sep=comma] {data_cube/FreqsNeumProblem_recaca.csv};
    %     \addplot [blue,mark=+] table [x=h, y=f_rhs_mean, col sep=comma] {data_cube/FreqsNeumProblem_recaca.csv};
    %     \addplot [blue!20] fill between[of=maxfreq and minfreq]; 
    %     % 
    %     \addplot [red] table [x=h, y=frhs_min, col sep=comma] {data_cube/FreqsNeumProblem_nolwr.csv}; 
    %     \addplot [red] table [x=h, y=frhs_max, col sep=comma] {data_cube/FreqsNeumProblem_nolwr.csv};
    %     \addplot [red,mark=+] table [x=h, y=f_rhs_mean, col sep=comma]
    %     {data_cube/FreqsNeumProblem_nolwr.csv};
    %   \end{loglogaxis}
    % \end{tikzpicture}
    \includegraphics{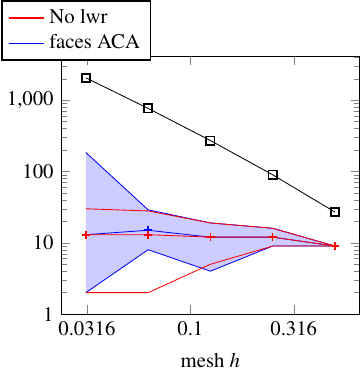}
  }
  \caption{\label{fig:CubeFreqsNeum}Cube: Number of used frequencies for the
    Neumann problem}
\end{figure}

After this check of a correct implementation and showing that the
compression does not spoil the results, the efficiency with respect to
the used frequencies, \ie the rank in the third matrix dimension is
studied. This gain of the 3D-ACA is presented in the
Figs.~\ref{fig:CubeFreqsDiri} and \ref{fig:CubeFreqsNeum} for the
Dirichlet and Neumann problem, respectively. In both figures, the
amount of used complex frequencies is given with respect to the mesh
size, \ie the refinement level.
Due to the selection of the different problems all of the four layer
potentials from above are studied. Note, as the 3D-ACA is applied on
the matrix blocks of the $\mathcal{H}$-matrix, each matrix block may
have its own number of  necessary complex frequencies. That
is why a min-, max- and mean-value is given.
Further, the numbers are displayed for the version using no
approximation in the faces (red) and those with an ACA approximation
within the faces (blue). Essentially, all operators show a drastic
reduce of the used frequencies resulting in a drastic reduction of
memory usage. Further, for more time steps (reducing time step size)
only a very moderate increase in the number of frequencies is
visible. Only for the double layer potentials the approximation with
ACA significantly changes the maximum number of used frequencies,
however not in the mean value. This is caused by the
structure of the double layer potential, which includes on the one
hand changing normal vectors resulting in blocks, which are badly
approximated by ACA, and
zero blocks. The zero blocks are not considered, however
some frequencies are tested to make sure that the block under
consideration is indeed a zero block. This results in the very low
number of the min-value. The higher max values result from
difficulties to determine the norm with \eqref{eq:norm} as the
monotonicity is not that good in the above discussed blocks.

\begin{figure}
  \centering
  \subcaptionbox{Single layer potential}{
    % \begin{tikzpicture}
    %   \begin{axis}[
    %     scale=.6, 
    %     colormap/cool,
    %     colorbar,
    %     xlabel=$Re(s_{\ell})$,
    %     ylabel= $Im(s_{\ell})$ ]
    %     \addplot [scatter,scatter src=explicit]
    %     table [x=Re,y=Im,meta=use, col sep=comma]
    %     {data_cube/slp_lv4.csv};
    %   \end{axis}
    % \end{tikzpicture}
    \includegraphics{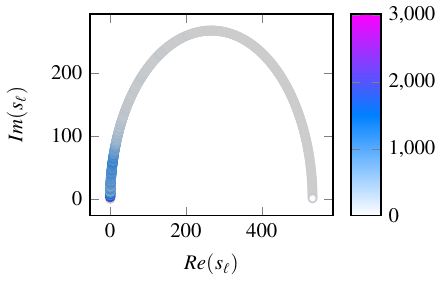}
  }\hfill
  \subcaptionbox{Hyper singular operator}{
    % \begin{tikzpicture}
    %   \begin{axis}[
    %     scale=.6, 
    %     colormap/cool,
    %     colorbar,
    %     xlabel=$Re(s_{\ell})$,
    %     ylabel= $Im(s_{\ell})$ ]
    %     \addplot [scatter,scatter src=explicit]
    %     table [x=Re,y=Im,meta=use, col sep=comma]
    %     {data_cube/hso_lv4.csv};
    %   \end{axis}
    % \end{tikzpicture}
    \includegraphics{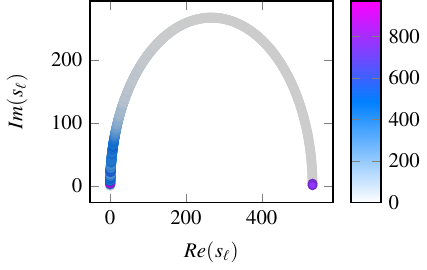}
  }
  \caption{\label{fig:CubeFreqsUsed}Cube (level 4): Used complex frequencies
    with a colour code for the number of matrix blocks, 
    where this frequency is active}
\end{figure}

To display the adaptively selected frequencies, in
Fig.~\ref{fig:CubeFreqsUsed} the complex frequencies of the gCQ are
plotted in the complex plane exemplarily for level 4.
A dense formulation would use all frequencies. The colour code shows
how much blocks in the 3D-ACA are evaluated at these frequencies. As
the algorithm starts for all blocks with the first frequency selected
by gCQ, this frequency is selected by all blocks (22432 for the single
layer potential and 8560 for the hyper sigular operator). This number is
deleted from the plot because else nothing would be visible due to the
linear scaled colour range. The figure clearly shows that there is a
region with small real parts where significant frequencies are
detected by the algorithm. Surprisingly, the last complex frequency
with large real part and nearly vanishing imaginary part is as well
selected in case of the hypersingular operator. The same happens for
the single layer potential, however for much less frequencies such
that it is not visible in the plots above. For both double layer
potentials similar pictures could be shown, where also for the adjoint
double layer potential the last frequency is more often selected. 

\begin{figure}
  \centering
  % \begin{tikzpicture}
  %   \begin{loglogaxis}
  %     [scale=.8,
  %     cycle list name=mylines,
  %     scale only axis,
  %     xlabel=mesh $h$,
  %     ylabel= compression,
  %     log ticks with fixed point,
  %     legend style={at={(0.7,.2)}}
  %     ]
  %     \addplot table [x=h, y=comp_part, col sep=comma] {data_cube/FreqsDiriProblem_nolwr.csv};
  %     \addlegendentry{Single layer potential: No lwr};
  %     \addplot table [x=h, y=comp_part, col sep=comma] {data_cube/FreqsDiriProblem_recaca.csv};
  %     \addlegendentry{Single layer potential: ACA};
  %     \addplot table [x=h, y=comp_part, col sep=comma] {data_cube/FreqsNeumProblem_nolwr.csv};
  %     \addlegendentry{Hypersingular operator: No lwr};
  %     \addplot table [x=h, y=comp_part, col sep=comma] {data_cube/FreqsNeumProblem_recaca.csv};
  %     \addlegendentry{Hypersingular operator: ACA};
  %   \end{loglogaxis}
  % \end{tikzpicture}
  \includegraphics{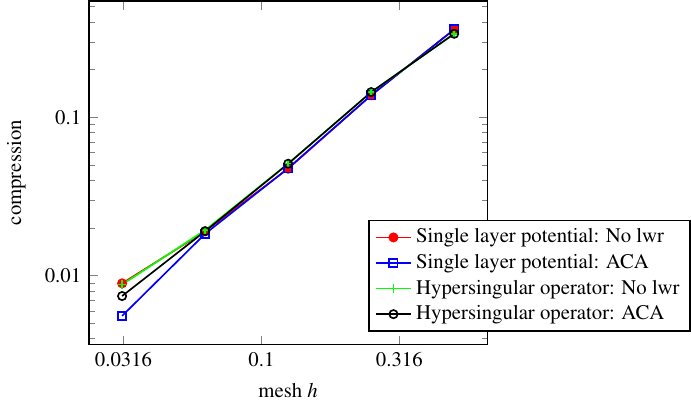}
  \caption{Cube: Compression of the data-sparse representation} 
  \label{fig:CubeCompress}
\end{figure}

The above results indicate that the memory requirement of the proposed
method is strongly improved compared to the standard formulation. As
the right hand side is multiplied on the fly, \ie the storage is the
same for the new and the original formulation, the compression is only
determined by the single layer and hypersingular operator in the
Dirichlet and Neumann problem, respectively.  This compression is
plotted in Fig.~\ref{fig:CubeCompress} for both problems, Dirichlet
and Neumann, versus the refinement. 
As mentioned above, the compression is the relation between the
required storage for the proposed formulation divided by the storage
of a dense formulation. Obviously, the improvement is very strong, where
the difference between the version without a compression in the face
to the results with the ACA applied in the faces is not that
strong. This is understandable if the sizes of the faces are
considered. Those are not that large such that ACA can achieve a
substantial improvement. Nevertheless, a compression of up to 0,56 \%
is a very good result.

\begin{figure}
  \centering
  % \begin{tikzpicture}
  %   \begin{loglogaxis}
  %     [scale=.8,
  %     cycle list name=mylines,
  %     scale only axis,
  %     xlabel=dofs,
  %     ylabel= CPU time,
  %     log ticks with fixed point,
  %     legend style={at={(.9,.35)}}
  %     ]
  %     \addplot table [x=nu, y=Diri_dense, col sep=comma] {data_cube/Timing.csv};
  %     \addlegendentry{Dirichlet: dense};
  %     \addplot table [x=nu, y=Diri_nolwr, col sep=comma] {data_cube/Timing.csv};
  %     \addlegendentry{Dirichlet: No lwr};
  %     \addplot table [x=nu, y=Diri_recaca, col sep=comma] {data_cube/Timing.csv};
  %     \addlegendentry{Dirichlet: ACA};
  %     \addplot table [x=du, y=Neum_dense, col sep=comma] {data_cube/Timing.csv};
  %     \addlegendentry{Neumann: dense};
  %     \addplot table [x=du, y=Neum_nolwr, col sep=comma] {data_cube/Timing.csv};
  %     \addlegendentry{Neumann: No lwr};
  %     \addplot table [x=du, y=Neum_recaca, col sep=comma] {data_cube/Timing.csv};
  %     \addlegendentry{Neumann: ACA};
  %     %
  %     \addplot [dashed,black] table [restrict expr to  domain={\thisrow{du}}{100:20000},x=du, y
  %     expr=5.e-2*\thisrow{N}*\thisrow{du}*ln(\thisrow{du})/ln(10), col sep=comma]
  %     {data_cube/FreqsNeumProblem_recaca.csv};
  %     \addlegendentry{complexity $\mathcal{O}(N M \log(M))$};
  %     \addplot [dotted,black] table [restrict expr to  domain={\thisrow{du}}{100:20000},x=du, y
  %     expr=3.e-6*\thisrow{N}*ln(\thisrow{N})/ln(10)*\thisrow{du}*\thisrow{du}, col sep=comma]
  %     {data_cube/FreqsNeumProblem_recaca.csv};
  %     \addlegendentry{complexity $\mathcal{O}(N (\log{N})^2 M^2)$};
  %   \end{loglogaxis}
  % \end{tikzpicture}
  \includegraphics{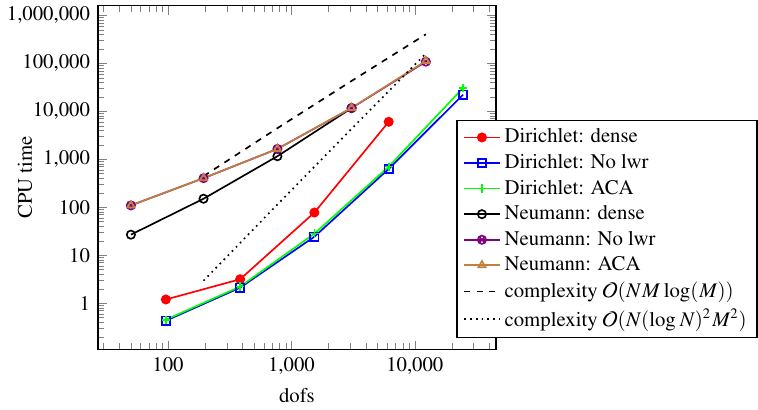}
  \caption{Cube: CPU time of the data-sparse representation}
  \label{fig:CubeTime}
\end{figure}

Concerning the computing time a not that clear picture is
possible. The presented results have been computed on one
  cluster node with 64 cores and 1 TB of storage. The code uses OpenMP
  for the obvious parallelisation with respect to the
  $\mathcal{H}$-matrix blocks. First, the complexity of the dense algorithm is
$\mathcal{O}(N (\log(N))^2 M^2)$. This line is plotted in
Fig.~\ref{fig:CubeTime} as dotted line. 
The dense Dirichlet problem follows this line, however, the Neumann
problem not that clear. The reason might be the Galerkin formulation,
which dominates the spatial matrix element calculation, and maybe
avoids to be in the asymptotic range. The data sparse computations a
clearly faster for the Dirichlet problem compared to the dense one but
for the Neumann problem the advantage is not yet visible. Again it
should be remarked that the Neumann problem uses a Galerkin method and
the Dirichlet problem a collocation method. For the 3D-ACA the
complexity line $\mathcal{O}(N M \log(M)$ seems to fit. However, it
must be remarked that there is not any reason why this should
hold. There is no indication in the algorithm that the adaptively
determined rank $r$ is linear in $N$. Overall, it can be concluded
that the 3D-ACA gives a very good performance with respect to the
storage reduction but not that much savings in the CPU-time.

\subsection{Example: Unit cube with reentrant corner (L-Shape) loaded by a smooth pulse}

\begin{figure}
  \subcaptionbox{L-Shape (level 1, $h=\unit[0.5]{m}$)}{
    \includegraphics[scale=.15,trim={20cm 2cm 20cm 7cm}, clip]{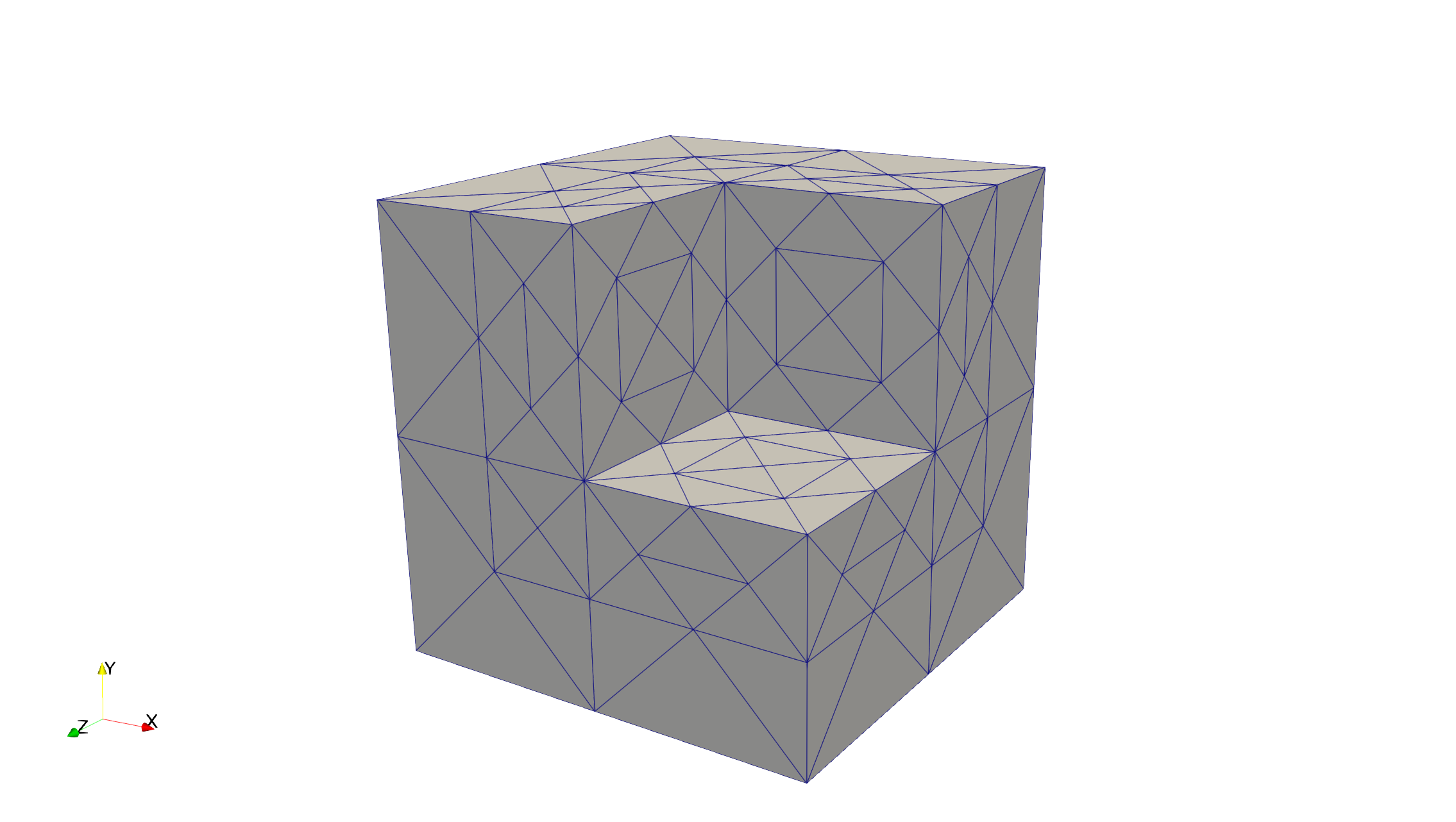}
  } \hfill
  \subcaptionbox{Used meshes}{
    \begin{tabular}{l|rrrr}
      level & nodes & elements & $h$ & $\Dt$\\ \hline
      1 & 95 & 186 & \unit[0.5]{m} & \unit[0.3]{s} \\
      2 & 374 & 744 & \unit[0.25]{m} & \unit[0.15]{s} \\
      3 & 1490 & 2976 & \unit[0.125]{m} & \unit[0.075]{s} \\
      4 & 5954 & 11904 & \unit[0.0625]{m} & \unit[0.0375]{s} \\
      5 & 23810 & 47616 & \unit[0.03125]{m} & \unit[0.01875]{s}
    \end{tabular}
  }
  \caption{\label{fig:LShape}L-Shape: Geometry and discretisation parameters}
\end{figure}

The next example is the same study as above but with a different
geometry. The selected geometry is the three-dimensional version of an
L-shape, \ie a cube with a reentrant corner. The geometry as well as
the parameters of the used meshes are given in Fig.~\ref{fig:LShape}.
The total observation time $T=\unit[3]{s}$ is selected such that the
smooth pulse travels over the whole computing domain. 
Note, the value $h$ in the table in Fig. ~\ref{fig:LShape} is
misleading as it gives values of the largest elements, which are located on
the not visible sides. Those elements on the side edges of the L-Shape
on the left and right side have these sizes as well.  Hence, towards the
reentrant corner the mesh is refined compared to $h$. As above, a Dirichlet and a
Neumann problem is studied. For both, the excitation point is selected
in the reentrant part at $\mathbf{x}=(0.25,  0.25, 0.25)^{\mathsf{T}}$
and the load is the same smooth puls of \eqref{eq:load}. Note, the
origin of the coordinate system is again in the center of the unit
cube. Different to
the setting above the Neumann problem is now not a scattering problem
but computes the pressure at the surface of this L-Shape. The used
equations are again \eqref{eq:bie_discr_diri} for the Dirichlet
problem solved with a collocation approach and the Neumann problem uses
\eqref{eq:bie_discr_neum} with a Galerkin approach.

\begin{figure}[hbt]
  \centering
  \includegraphics{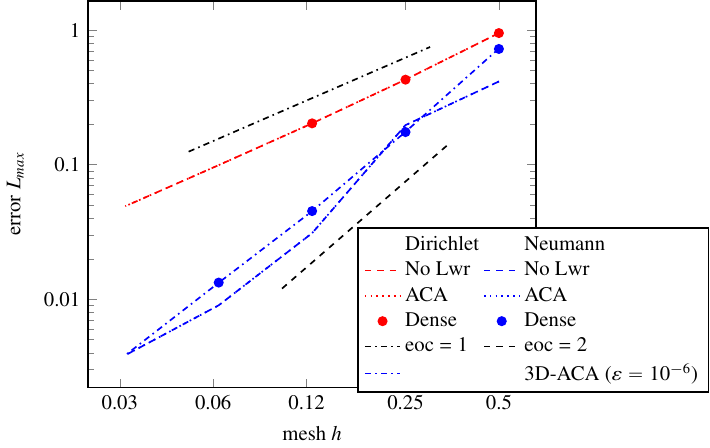}
  \caption{L-Shape: $L_{max}$-error versus refinement in space and time}
  \label{fig:ConvLShape}
\end{figure}

First, as in the test above the convergence rates are reported to show
that the proposed method keeps the error level of the dense
computation. In Fig.~\ref{fig:ConvLShape}, the $L_{max}$-error is
plotted versus the refinement in space and time for the Dirichlet and
Neumann problem.
In this figure three lines are visible. The Dirichlet problem is
solved with the expected convergence rate and the dense computation
gives the same results. For the Neumann problem the rate is
still obtained also with the 3D-ACA but the error level is
different. However, as indicated with the line '3D-ACA
($\varepsilon=10^{-6}$)' it is shown that for an increased precession
the dense results can be recovered. This result confirms that the 3D-ACA
has an influence on the computation and the $\varepsilon$-parameter must
be carefully adjusted, despite that in this case the coarser
approximation gives better results.

\begin{figure}
  \subcaptionbox{Single layer potential}{
    % \begin{tikzpicture}
    %   \begin{loglogaxis}
    %     [scale=.6,
    %     scale only axis,
    %     xlabel=mesh $h$,
    %     ylabel= Number frequencies,
    %     log ticks with fixed point,
    %     legend columns=2,
    %     legend style={at={(0,1.1)}}
    %     ]
    %     \addplot [black,mark=square] table [x=h, y=no_f, col sep=comma] {data_Lshape/FreqsDiriProblem_nolwr.csv}; 
    %     \addlegendentry{No 3D-ACA};
    %     % 
    %     \addplot [name path = minfreq,blue] table [x=h, y=fmat_min, col sep=comma] {data_Lshape/FreqsDiriProblem_recaca.csv};
    %     \addlegendentry{Min no Freqs};
    %     \addplot [name path = maxfreq, blue] table [x=h, y=fmat_max, col sep=comma] {data_Lshape/FreqsDiriProblem_recaca.csv};
    %     \addlegendentry{Max no Freqs};
    %     \addplot [blue,mark=+] table [x=h, y=fmat_mean, col sep=comma] {data_Lshape/FreqsDiriProblem_recaca.csv};
    %     \addlegendentry{Mean no Freqs};
    %     \addplot [blue!20] fill between[of=maxfreq and minfreq]; 
    %     % 
    %     \addplot [red] table [x=h, y=fmat_min, col sep=comma] {data_Lshape/FreqsDiriProblem_nolwr.csv}; 
    %     \addplot [red] table [x=h, y=fmat_max, col sep=comma] {data_Lshape/FreqsDiriProblem_nolwr.csv};
    %     \addplot [red,mark=+] table [x=h, y=fmat_mean, col sep=comma]
    %     {data_Lshape/FreqsDiriProblem_nolwr.csv};
    %   \end{loglogaxis}
    % \end{tikzpicture}
    \includegraphics{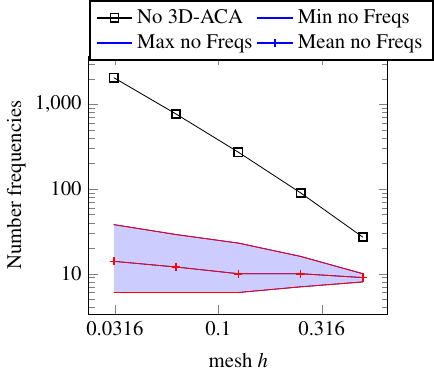}
  }
  \subcaptionbox{Double layer potential}{
    % \begin{tikzpicture}
    %   \begin{loglogaxis}
    %     [scale=.6, 
    %     scale only axis,
    %     xlabel=mesh $h$,
    %     log ticks with fixed point,
    %     legend style={at={(-.2,1.1)}}
    %     ]
    %     \addlegendimage{red} \addlegendentry{No lwr};
    %     \addlegendimage{blue} \addlegendentry{faces ACA};
    %     \addplot [black,mark=square] table [x=h, y=no_f, col sep=comma] {data_Lshape/FreqsDiriProblem_nolwr.csv}; 
    %     % 
    %     \addplot [name path = minfreq,blue] table [x=h, y=frhs_min, col sep=comma] {data_Lshape/FreqsDiriProblem_recaca.csv};
    %     \addplot [name path = maxfreq, blue] table [x=h, y=frhs_max, col sep=comma] {data_Lshape/FreqsDiriProblem_recaca.csv};
    %     \addplot [blue,mark=+] table [x=h, y=f_rhs_mean, col sep=comma] {data_Lshape/FreqsDiriProblem_recaca.csv};
    %     \addplot [blue!20] fill between[of=maxfreq and minfreq]; 
    %     % 
    %     \addplot [red] table [x=h, y=frhs_min, col sep=comma] {data_Lshape/FreqsDiriProblem_nolwr.csv}; 
    %     \addplot [red] table [x=h, y=frhs_max, col sep=comma] {data_Lshape/FreqsDiriProblem_nolwr.csv};
    %     \addplot [red,mark=+] table [x=h, y=f_rhs_mean, col sep=comma]
    %     {data_Lshape/FreqsDiriProblem_nolwr.csv};
    %   \end{loglogaxis}
    % \end{tikzpicture}
    \includegraphics{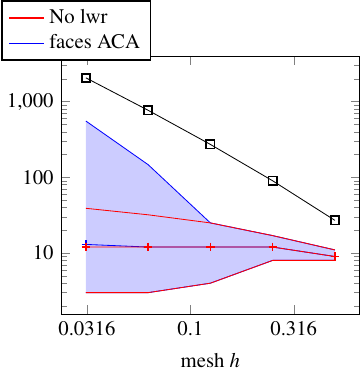}
  }
  \caption{\label{fig:LShapeFreqsDiri}L-Shape: Number of used frequencies for the
    Dirichlet problem}
\end{figure}

\begin{figure}
  \subcaptionbox{Hyper singular operator}{
    \includegraphics{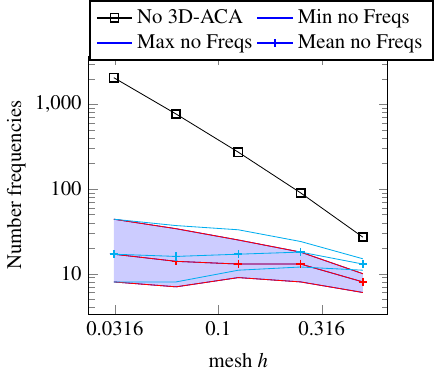}
  }
  \subcaptionbox{Adjoint double layer potential}{
    \includegraphics{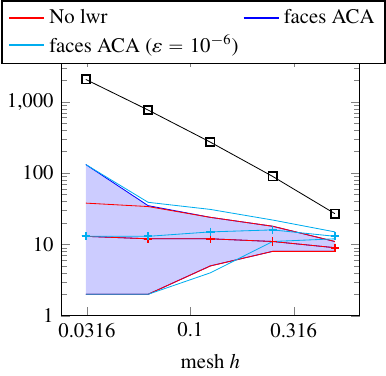}
  }
  \caption{\label{fig:LShapeFreqsNeum}L-Shape: Number of used frequencies for the
    Neumann problem}
\end{figure}

The used complex frequencies are as well studied for this
geometry. They are displayed in Figs.~\ref{fig:LShapeFreqsDiri} and \ref{fig:LShapeFreqsNeum}. 
Essentially, the pictures are similar to those for the unit
cube and the numbers are comparable. Hence, there is no evidence that
the geometry affects the used frequencies. 
Another remarkable point is the different numbers for an increased
precession value $\varepsilon$. For the smaller meshes and less time
steps the difference in the mean value is clearly visible but with
increasing numbers for the times steps this difference
vanishes. Hence, it may be concluded that the adaptivity inherently in
the proposed algorithm detects the amount of frequencies necessary to
give the physical behaviour of the solution. The compression rates are
comparable to the unit cube and also the selection of the complex
frequencies. Hence, these figures are skipped for the sake of brevity.

\subsection{Example: Unit ball (sphere) loaded by a non-smooth pulse}
The next example is used to show the influence of the time step size
on the proposed method. The main feature of gCQ is to allow variable
time step sizes to improve the approximation of solutions. Such
gradings in time are necessary for non-smooth time behaviour of the
solution. In the paper studying the numerical behaviour of the
gCQ~\cite{lopezfernandez12a}, a solution is selected where the behaviour
in time is non-smooth such that a grading in the time mesh makes
sense. This example is also used here to show the ability of the
proposed approximation method to handle as well such solutions.

\begin{figure}[hbt]
  \centering
  \subcaptionbox{Sphere with radius $\unit[1]{m}$}{
    \includegraphics[scale=.25]{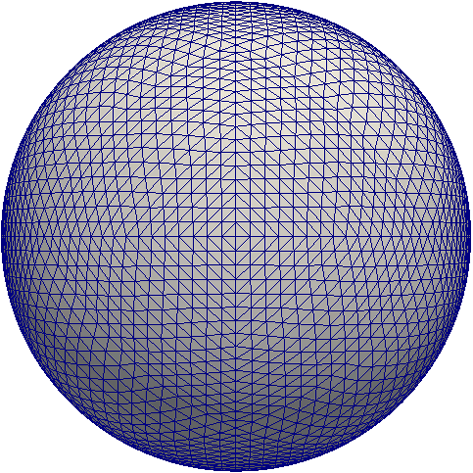}
  } \hspace{2cm}
  \subcaptionbox{Used mesh and time step sizes}{
    \begin{tabular}{l|rrrr}
      level & elements & $h$ & $\Dt_{const}$ \\ \hline
      1 & 7648 & $\approx \unit[0.05]{m}$ & \unit[0.15]{s} \\
      2 & 7648 & $\approx \unit[0.05]{m}$  & \unit[0.075]{s} \\
      3 & 7648 & $\approx \unit[0.05]{m}$  & \unit[0.0375]{s} \\
      4 & 7648 & $\approx \unit[0.05]{m}$  & \unit[0.01875]{s} \\
      5 & 7648 & $\approx \unit[0.05]{m}$  & \unit[0.00938]{s}
    \end{tabular}
  }
  \caption{\label{fig:sphere}Sphere: Geometry, mesh size and used time
    step sizes}
\end{figure}

\begin{figure}
  \centering
  % \begin{tikzpicture}
  %   \begin{loglogaxis}
  %     [scale=.75,
  %     cycle list name=mylines,
  %     scale only axis,
  %     xlabel=time step $\Dt_{const}$,
  %     ylabel= error $L_{max}$,
  %     log ticks with fixed point,
  %     legend columns=2,
  %     legend style={at={(.7,.17)}}
  %     ]
  %     \addlegendimage{empty legend};
  %     \addlegendentry{constant $\alpha=1$};
  %     \addlegendimage{empty legend};\addlegendentry{variable $\alpha=2$}; 
  %     %
  %     \addplot[blue, only marks,mark=square] table [x=Dt, y=Lmax_dense_a1, col sep=comma] {data_Sphere/ErrorsDiriProblem.csv};
  %     \addlegendentry{Dense};
  %     \addplot[red, only marks,mark=square] table [x=Dt, y=Lmax_dense_a2, col sep=comma] {data_Sphere/ErrorsDiriProblem.csv};
  %     \addlegendentry{Dense};
  %     % 
  %     \addplot[blue,mark=|] table [x=Dt, y=Lmax_recaca_a1, col sep=comma] {data_Sphere/ErrorsDiriProblem.csv};
  %     \addlegendentry{3D-ACA};
  %     \addplot[red,mark=|] table [x=Dt, y=Lmax_recaca_a2, col sep=comma] {data_Sphere/ErrorsDiriProblem.csv};
  %     \addlegendentry{3D-ACA};
  %     % 
  %     \addplot[dashdotted] gnuplot[domain=.02:.1] {2*x**(0.5)};
  %     \addlegendentry{eoc = 0.5};
  %     \addplot[dashed] gnuplot[domain=.03:.1] {3*x};
  %     \addlegendentry{eoc = 1};
  %   \end{loglogaxis}
  % \end{tikzpicture}
  \includegraphics{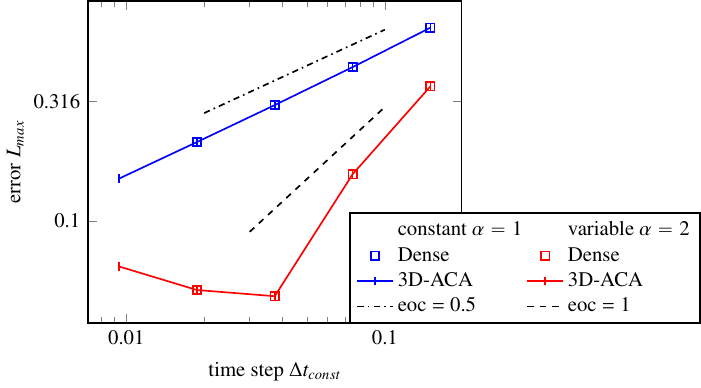}
  \caption{Sphere: $L_{max}$-error versus time step size for constant
  and variable time step size}
  \label{fig:ConvSphere}
\end{figure}

The selected geometry is a ball (sphere) with radius $\unit[1]{m}$. The
geometry and the selected mesh are shown in Fig.~\ref{fig:sphere}.
Different to the examples above, here, the spatial
discretisation is 
not changed. This is motivated by the constructed solution. It is
known that the density function within a single layer approach for the
sphere corresponds to the spherical harmonics times the time
derivative of the applied excitation in
time~\cite{sauterveit13b}. Selecting the zeroth order spherical
harmonic $Y_0^0\kl{\x}$, the spatial solution is a constant
function. This part of the solution is 
decoupled from the temporal solution, which might justify to
keep the spatial discretisation constant and refine only in time. The
selected excitation and analytical solution is
\begin{equation} \label{eq:load_sphere}
  g_D \kl{\x,t} = Y_0^0\kl{\x} t^{\frac{3}{2}} e^{-t} \quad
  \Rightarrow \quad \vartheta\kl{\x,t}  = Y_0^0\kl{\x}
  \sum\limits_{k=0}^{\lfloor t/2 \rfloor} e^{-(t-2k)} \kl{\frac{3}{2}
    \kl{t-2k}^{\frac{1}{2}} - \kl{t-2k}^{\frac{3}{2} }} \; .
\end{equation}
The term with the square root shows a non-smooth behaviour in time
starting at $t=\unit[0]{s}$ and is repeated in a periodic
manner. Selecting the total observation time $T=\unit[1.5]{s}$ results
in a solution, where the non-smooth part is only at the
beginning. Hence, a graded time mesh
\begin{equation} \label{eq:timeMesh}
  t_n = T \kl{\frac{n}{N}}^{\alpha} = N \Dt_{const}
  \kl{\frac{n}{N}}^{\alpha} \qquad n=0, 1, \ldots , N
\end{equation}
can be recommended. The value $\Dt_{const}$ is introduced to clarify
the meaning of time step size, \ie in Fig.~\ref{fig:sphere} this value
is listed. For $\alpha=2$ the full convergence order should
be preserved. The effect of the grading can be observed in
Fig.~\ref{fig:ConvSphere}, where the $L_{max}$-error is presented
versus a refinement in time for different exponents $\alpha$.
The value $\alpha=1$ corresponds to a constant step size, where as
expected the convergence order decreases to $\mathrm{eoc}=0.5$. With
the graded mesh the convergence order is recovered. Note, the density
function is approximated by constant shape functions, hence a order
larger than 1 can not be expected. Obviously, there is a break in the
convergence, where a further decrease of the time step size does not
improve the result. This behaviour is already reported
in~\cite{lopezfernandez12a} and reflects that the spatial
discretisation is no longer suitable in combination with the spatial
integration. Using the analytical solution in the spatial variable
allows gCQ to keep the convergence order
(see~\cite{lopezfernandez12a}). The proposed accelerated algorithm
shows obviously the same behaviour. Hence, the approximation due to the
3D-ACA does not affect the efficiency of the gCQ based time
domain BEM.

\begin{figure}[hbt]
  \centering
  % \begin{tikzpicture}
  %   \begin{loglogaxis}
  %     [scale=.8,
  %     cycle list name=mylines,
  %     scale only axis,
  %     xlabel=time step size $\Dt$,
  %     ylabel= Number frequencies,
  %     log ticks with fixed point,
  %     legend style={at={(0.5,.85)}}
  %     ]
  %     \addplot [black,mark=square] table [x=Dt, y=no_f, col sep=comma] {data_Sphere/FreqsDiriProblem.csv}; 
  %     \addlegendentry{No 3D-ACA};
  %     % 
  %     \addplot [name path = minfreq,red] table [x=Dt, y=al2_min, col sep=comma] {data_Sphere/FreqsDiriProblem.csv}; 
  %     \addlegendentry{Min no Freqs (red $\alpha=2$, blue $\alpha=1$)};
  %     \addplot [name path = maxfreq, red] table [x=Dt, y=al2_max, col sep=comma] {data_Sphere/FreqsDiriProblem.csv}; 
  %     \addlegendentry{Max no Freqs (red $\alpha=2$, blue $\alpha=1$)};
  %     \addplot [red,mark=+] table [x=Dt, y=al2_mean, col sep=comma] {data_Sphere/FreqsDiriProblem.csv}; 
  %     \addlegendentry{Mean no Freqs (red $\alpha=2$, blue $\alpha=1$)};
  %     \addplot [red!20] fill between[of=maxfreq and minfreq];
  %     %
  %     \addplot [name path = min2,blue] table [x=Dt, y=al1_min, col sep=comma] {data_Sphere/FreqsDiriProblem.csv};
  %     \addplot [name path = max2,blue] table [x=Dt, y=al1_max, col sep=comma] {data_Sphere/FreqsDiriProblem.csv};
  %     \addplot [blue,mark=+] table [x=Dt, y=al1_mean, col sep=comma]{data_Sphere/FreqsDiriProblem.csv};
  %     \addplot [blue!15, nearly opaque] fill between[of=max2 and min2];
  %   \end{loglogaxis}
  % \end{tikzpicture}
  \includegraphics{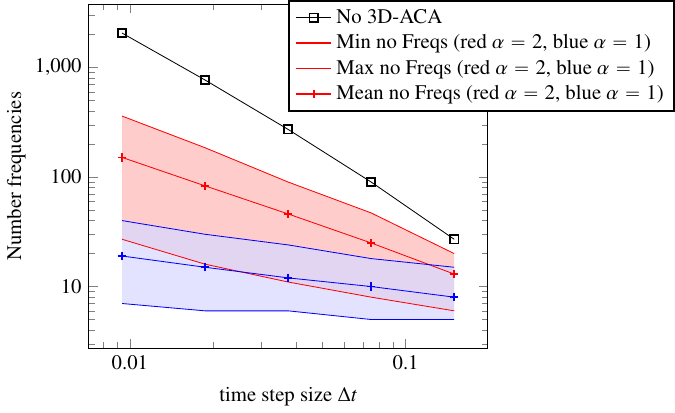}
  \caption{Sphere: Necessary frequencies for 3D-ACA where also the faces uses ACA}
  \label{fig:FreqSphere}
\end{figure}

Next, the used complex frequencies are studied for this example. In
Fig.~\ref{fig:FreqSphere}, the amount of used complex frequencies is
plotted for both temporal discretisations.
Obviously, the graded time mesh requires more frequencies compared to
the constant discretisation. This holds true in the mean-value as well
as in the max- and min-values. The tendency is as well that in this
case the mean value does not approach a nearly constant amount of time
steps. However, the increase is significantly smaler than in the dense
method and very good compression can be obtained. The compression is
for the first level 26\% and in the fifth level 0,98 \%. 
This reduction in storage is significant and allows to treat larger
problems. Unfortunately, for
  $\alpha=2$ the computation time increases compared to the uniform time step. If the same
  accuracy is observed, \ie comparing the computing time of refinement
  level 2 for $\alpha=2$ with level 5 for $\alpha=1$ (see
  Fig.~\ref{fig:ConvSphere}), almost a factor of 2 is measured
  (\unit[7746]{sec} to \unit[4295]{sec}). The variable step
  size does not pay off in this example even when using 3D-ACA.

\begin{figure}
  \centering
  \subcaptionbox{$\alpha=1$}{
    % \begin{tikzpicture}
    %   \begin{axis}[
    %     scale=.63,
    %     colormap/cool,
    %     colorbar,
    %     xlabel=$Re(s_{\ell})$,
    %     ylabel= $Im(s_{\ell})$ ]
    %     \addplot [scatter,scatter src=explicit]
    %     table [x=Re,y=Im,meta=use, col sep=comma]
    %     {data_Sphere/alpha1.csv};
    %   \end{axis}
    % \end{tikzpicture}
    \includegraphics{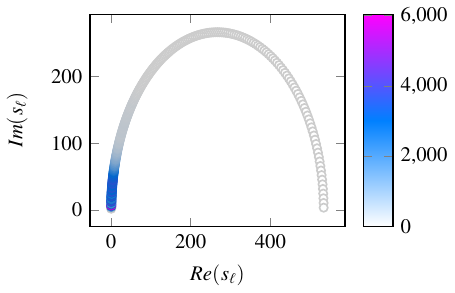}
  }
  \subcaptionbox{$\alpha=2$}{
    % \begin{tikzpicture}
    %   \begin{axis}[
    %     scale=.63,
    %     colormap/cool,
    %     colorbar,
    %     xlabel=$Re(s_{\ell})$,
    %     ylabel= $Im(s_{\ell})$ ]
    %     \addplot [scatter,scatter src=explicit]
    %     table [x=Re,y=Im,meta=use, col sep=comma]
    %     {data_Sphere/alpha2.csv};
    %   \end{axis}
    % \end{tikzpicture}
    \includegraphics{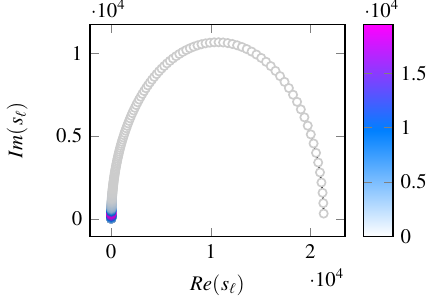}
  }
  \caption{\label{fig:sphereColorFreqs}Sphere (level 3): Used complex frequencies
    with a colour code for the number of matrix blocks, 
    where this frequency is active }
\end{figure}

For this example it is interesting to report not only the
numbers of used frequencies but also which frequencies are
selected. In Fig.~\ref{fig:sphereColorFreqs}, the same plots as
presented for the unit cube are shown.
First, looking at the right plot for the graded time mesh $\alpha=2$ it can be
clearly seen how focused the frequencies are towards the imaginary
axis and small real values. It is even visible that for large real
values a coarsening happens. This is in accordance with the
integration rule. The difference between the largest and smallest time
step (parameter $q$ in~\ref{app:gcq}) determines the
concentration of the frequencies towards the imaginary axis.  Another
effect independent of the 3D-ACA is the much larger values of the
complex frequencies, \ie the circle is larger. Please, have a look at
the values on the axis compared to the left picture. Somehow as a
consequence of this concentration as well the selected frequencies are
concentrated towards small real values.  It seems that the 3D-ACA
preserves the properties of the integration rule.
\section{Conclusions}
Time domain boundary element formulations are usually very storage
demanding as a lot of matrices for several time steps have to be
stored. This is as well a significant drawback for the generalised convolution
quadrature (gCQ) based formulation, where approximately $N \log N$ or
more matrices each of the size of an elliptic problem have to be
stored or more precise kept in memory. This reduces possible
applications to very small sizes. The generalised adaptive cross
approximation (3D-ACA) allows to find a low rank representations of
three-dimensional tensors, hence this method is applied here to a gCQ based
time domain BEM.

The results show a very significant reduction in storage while keeping
the convergence of a dense computation. The algorithm is somehow a
black-box technique as no modifications of the existing kernels of the dense
formulation is necessary. The presented studies show the storage
reduction for three different geometries, a unit cube, an L-shape, and
a sphere. All three show similar behaviour with respect to the
adaptively selected complex frequencies of the gCQ. The example with
the sphere shows as well that the benefit of gCQ to select non-uniform
time steps can be preserved, however for strongly different time step
sizes the amount of necessary frequencies increases. Summarising, the
application of 3D-ACA on the gCQ based time domain BEM gives a data
sparse method with storage savings of more than 99\%.

\paragraph*{Acknowledgement}
The work of M. Schanz is partially supported by the joint DFG/FWF Collaborative Research Centre CREATOR (DFG: Project-ID 492661287/TRR 361; FWF: 10.55776/F90) at TU Darmstadt, TU Graz and JKU Linz.

\appendix
\section{Parameters of the gCQ} \label{app:gcq}
% Parameters for gCQ
The derivation and reasoning how the integration weights and points
are determined can be found in~\cite{lopezfernandez12a, lopezfernandez15b}. The result of
these papers are recalled here.
The integration points in the complex plane are
\begin{align*}
  s_\ell = \gamma(\sigma_\ell)\ , & & \omega_\ell = \frac{4 K(k^2)}{1
                                      \pi i} \gamma^\prime
                                      (\sigma_\ell)\ , && N_Q = N \log(N)\ , 
\end{align*}
where for Runge-Kutta methods with  $m>1$  stages it should be $N_Q = N
(\log(N))^2$. $K(k)$ is the complete elliptic integral of first kind
\begin{align*}
 K(k) = \int^1_0 \frac{dx}{\sqrt{(1-x^2) (1-k^2 x^2)}}\ , & & K^\prime (k) = K(1-k)
\end{align*}
and $K^\prime$ is its derivative, which equals the integral of the complementary modulus.
The argument $k$ depends on the relation $q$ of the maximum and minimum step sizes in the following way
\begin{align*}
 k = \frac{q - \sqrt{2q-1}}{q + \sqrt{2q-1}}  \qquad q = 5
  \frac{\Dt_{\text{max}} \max_i
  |\lambda_i\kl{\A}|}{\Dt_{\text{min}}\min_i
  |\lambda_i\kl{\A}|} \; ,
\end{align*}
with the eigenvalues $\lambda_i\kl{\A}$. For the implicit Euler method
the eigenvalues are 1 and the factor 5 in $q$ can be skipped. 
The functions $\gamma(\sigma_\ell)$ and $\gamma^\prime(\sigma_\ell)$ are
\begin{align*}
 \gamma(\sigma_\ell) =& \frac{1}{\Dt_{\text{min}}(q-1)} \left(
                        \sqrt{2q-1} \frac{k^{-1} +
                        \text{sn}(\sigma_\ell)}{k^{-1}-\text{sn}(\sigma_\ell)}
                        -1 \right)  \\[1ex]
 \gamma^\prime(\sigma_\ell) =&
                               \frac{\sqrt{2q-1}}{\Dt_{\text{min}}(q-1)}
                               \frac{2\ \text{cn}(\sigma_\ell)\
                               \text{dn}(\sigma_\ell)}{k (k^{-1}-
                               \text{sn}(\sigma_\ell))^2}  \\[1ex]
 \sigma_\ell =& -K(k^2) + \left(\ell-\frac{1}{2} \right) \frac{4
                K(k^2)}{N_Q} + \frac{i}{2} K^\prime(k^2) \, ,
\end{align*}
where $\text{sn}(\sigma_\ell)$, $\text{cn}(\sigma_\ell)$ and
$\text{dn}(\sigma_\ell)$ are the Jacobi elliptic functions. 
As seen from above, the integration contour is only determined by the
largest and the smallest time steps chosen but not dependent on any
intermediate step sizes. Due to the symmetric distribution of the integration points with
respect to the real axis, only half of the frequencies $s_\ell$ need to be
calculated.

Last, it may be remarked that for constant time steps $\Dt_i=const$ the parameter determination would fail
because this choice results in $q=1$ and $k=0$. Unfortunately, this
value is not allowed for the complete elliptic integral. However, a
slight change in the parameter $q \approx 1$ fixes this problem
without spoiling the algorithm. The latter can be done as these
parameter choices are empirical.

%\bibliography{string_abbrev, scha_Journal, scha_Conference, scha_Book, fastBEM,bem,allgemein}
% \bibliographystyle{ieeetr}
%\bibliographystyle{plainnat}

\end{document}